\documentclass[12pt, a4paper]{article}
\usepackage[left=1in,right=1in,top=1.0in,bottom=1in]{geometry}
\usepackage[T1]{fontenc}
\usepackage[utf8]{inputenc}
\usepackage{color}
\usepackage{float}
\usepackage{amsmath}
\usepackage{amsthm}
\usepackage{amssymb}
\usepackage{graphicx}
\usepackage{esint}

\makeatletter
\numberwithin{equation}{section}
\numberwithin{figure}{section}

\usepackage{caption}
\usepackage{subcaption}
\usepackage{authblk}
\usepackage{subfig}

\usepackage[numbers]{natbib}
\usepackage{authblk} 

\makeatother

\usepackage{babel}
\begin{document}
\title{A finite element model for thermomechanical stress-strain fields in transversely isotropic strain-limiting materials}

\author[1]{Saugata Ghosh}
\author[2]{Dambaru Bhatta}
\author[3]{S. M. Mallikarjunaiah} 

\affil[1,2]{School of Mathematical \& Statistical Sciences,\\
The University of Texas - Rio Grande Valley, \\
Edinburg, Texas 78539, USA \\
Email:  saugata.ghosh01@utrgv.edu, dambaru.bhatta@utrgv.edu}

\affil[3]{Department of Mathematics \& Statistics,\\
Texas A\&M University-Corpus Christi, \\
Corpus Christi, Texas 78412-5825, USA \\
Email: m.muddamallappa@tamucc.edu (Corresponding Author)}
\date{}

\maketitle
\begin{abstract}
This paper presents a comprehensive computational framework for investigating thermo-elastic fracture in transversely isotropic materials, where classical linear elasticity fails to predict physically realistic behavior near stress concentrations. We address the challenge of unphysical strain singularities at crack tips by employing a strain-limiting theory of elasticity. This theory is characterized by an algebraically nonlinear constitutive relationship between stress and strain, which intrinsically enforces a limit on the norm of the strain tensor. This approach allows the development of very large stresses, as expected near a crack tip, while ensuring that the corresponding strains remain physically bounded. A loosely coupled system of linear and nonlinear partial differential equations governing the response of a thermo-mechanical transversely isotropic solid is formulated. We develop a robust numerical solution based on the finite element method, utilizing a conforming finite element discretization within a continuous Galerkin framework to solve the two-dimensional boundary value problem. The model is applied to analyze the stress and strain fields near an edge crack under severe thermo-mechanical loading. Our numerical results reveal a significant departure from classical predictions: while stress concentrates intensely at the crack tip, the strain grows at a substantially slower rate and remains bounded throughout the domain. This work validates the efficacy of the strain-limiting model in regularizing thermo-elastic crack-tip fields and establishes a reliable computational foundation for the predictive modeling of thermally driven crack initiation and evolution in advanced anisotropic materials.
\end{abstract}

\section{Introduction}

The confluence of extreme thermal and mechanical loads presents a formidable challenge in modern engineering, particularly in the design and analysis of high-performance systems. In fields ranging from aerospace to energy production, structural components are routinely subjected to environments where intense heat and significant mechanical stress interact. A critical aspect of this challenge lies in the material's constitutive response, which is often not uniform in all directions. Many advanced materials, such as fiber-reinforced composites, laminated ceramics, and certain geological formations, exhibit \textit{transverse isotropy}, where their mechanical and thermal properties are directionally dependent. For instance, jet engine turbine blades must withstand extreme temperatures while enduring immense centrifugal forces, and thermal barrier coatings on their surfaces are designed to resist thermal shock that can lead to delamination and fracture \cite{Clarke2003a,Clarke2003b,krishnasamy2018}. Similarly, spacecraft re-entering the atmosphere experience a combination of severe aerodynamic pressure and rapid thermal gradients.  In these applications, the material's anisotropic nature is not a minor detail but a dominant factor that governs deformation, damage initiation, and potential failure. Therefore, developing predictive models that can accurately capture stress and strain fields in such multi-physical, anisotropic settings is essential for ensuring the safety, reliability, and longevity of these critical structures.

The classical theories of thermoelasticity and linear elastic fracture mechanics (LEFM) have long served as the foundational frameworks for analyzing such problems \cite{Anderson2005, Inglis1913, love1944treatise}. However, they suffer from fundamental limitations that manifest as non-physical artifacts, especially under extreme conditions. First, classical thermoelasticity, which is based on Fourier's law of heat conduction, assumes that thermal disturbances propagate at an infinite speed---a mathematical convenience that breaks down in applications involving rapid thermal transients or high-frequency phenomena. To resolve this, generalized thermoelastic theories were developed, such as those proposed by Lord and Shulman \cite{lord1967} and Green and Lindsay \cite{Green1972}. These models, which build upon the work of Cattaneo \cite{Cattaneo1958} and Vernotte \cite{Vernotte1958}, introduce a thermal relaxation time to ensure a finite speed for heat waves. While these advanced theories provide a more realistic description of thermal dynamics, they inherit a second major flaw from classical mechanics: the prediction of infinite strain at geometric singularities like crack tips \cite{broberg1999}. This mathematical singularity is physically untenable and complicates the prediction of fracture initiation and propagation.

To address the challenge of unbounded crack-tip strains, a more sophisticated class of constitutive models is required. A significant advancement in this area is the \textit{implicit constitutive theory} introduced by Rajagopal \cite{rajagopal2003implicit}, which offers a powerful generalization of classical elasticity. Within this framework, the stress and strain tensors are related through an implicit equation, enabling the description of a much broader range of material behaviors beyond simple linear elasticity. A particularly important subclass of these models is known as \textit{strain-limiting elasticity} \cite{rajagopal2007elasticity,rajagopal2011modeling,gou2015modeling}. These models are designed to ensure that strain remains finite throughout the material, even in the presence of infinite stress concentrations. By positing strain as a nonlinear function of stress, they effectively introduce a constitutive ``ceiling'' that regularizes the strain field near discontinuities \cite{bulivcek2014elastic, bustamante2009some, itou2018states, kulvait2013, Mallikarjunaiah2015}. This feature aligns the mathematical prediction with physical reality, where materials cannot deform infinitely.

Recent research has successfully demonstrated the efficacy of this approach by developing a finite element framework for coupled thermo-mechanical problems in \textit{isotropic} bodies, showing that the strain-limiting model effectively regularizes strain growth near crack tips \cite{yoon2022finite}. However, as previously noted, many materials used in demanding applications are not isotropic but are instead characterized by {transverse isotropy} \cite{kachanov1993}. Applying an isotropic model to such materials is an oversimplification that can mask the critical directional effects governing the onset and evolution of failure. The extension of strain-limiting thermoelastic models to transversely isotropic media is therefore a natural and necessary step toward developing more accurate and reliable engineering tools.

This work bridges this critical gap by formulating and analyzing a finite element model for a coupled thermo-mechanical system within a transversely isotropic, strain-limiting thermoelastic body. Our contributions build directly upon the isotropic model introduced in \cite{yoon2022finite}, extending its theoretical and computational framework to account for directional dependencies in material behavior. The principal objective is to establish a robust rationale for employing nonlinear strain-limiting models in the context of transversely isotropic materials and to demonstrate their effectiveness in regulating strain growth near singularities. We analyze a two-dimensional system featuring one-way coupling between a linear steady-state heat conduction equation and a quasi-linear momentum balance equation. The resulting system is solved numerically using a finite element approach, where the nonlinearity arising from the constitutive law is handled with Picard's method. Through simulations of edge-cracked domains, we show that the nonlinear model produces a markedly slower growth of strain near the crack tip compared to stress, providing a more physically sound representation of material behavior in the presence of flaws.

The remainder of this paper is organized as follows. Section 2 establishes the necessary mathematical preliminaries and notational conventions. In Section 3, we introduce the theoretical framework of implicit constitutive relations, which is central to the development of the strain-limiting model. Building upon this foundation, Section 4 formulates the governing equations for the coupled thermo-elastic problem in transversely isotropic media. The numerical implementation, based on the Galerkin finite element method and a Picard iterative scheme for handling the nonlinearity, is detailed in Section 5. Section 6 presents and discusses the numerical results, including an analysis of how various model parameters influence the outcome. Finally, Section 7 offers concluding remarks and a summary of our key findings.

\section{Mathematical Preliminaries}

This section establishes the kinematic framework and notational conventions for analyzing a thermally conductive, transversely isotropic elastic body. The objective is to precisely define the measures of deformation and strain that are essential for modeling the material's thermo-mechanical response.

We consider a body that occupies an open, bounded, and connected domain $\Omega \subset \mathbb{R}^{2}$ in its reference (undeformed) configuration. The boundary of this domain is assumed to be Lipschitz continuous, which ensures it is sufficiently regular for the application of standard integral theorems and the formulation of well-posed boundary value problems. The vector space of second-order symmetric tensors in two dimensions is denoted by $Sym(\mathbb{R}^{2\times2})$. This space is equipped with the standard Frobenius inner product, defined for any two tensors $\mathbf{A}_{1}, \mathbf{A}_{2} \in Sym(\mathbb{R}^{2\times2})$ as $\mathbf{A}_{1}:\mathbf{A}_{2} = \mathrm{tr}(\mathbf{A}_{1}^{T}\mathbf{A}_{2})$, where $\mathrm{tr}(\cdot)$ is the trace operator. The associated norm is given by $\mathbf{\left\Vert A\right\Vert = \sqrt{A:A}}$.

The motion of the body is described by mapping each material point $\boldsymbol{X}$ from the reference configuration $\Omega$ to its corresponding spatial position $\boldsymbol{x}$ in the current (deformed) configuration. The {displacement vector} $\boldsymbol{u}:\Omega \longrightarrow \mathbb{R}^{2}$ quantifies the change in position of a material point and is defined as:
$$
\boldsymbol{u}(\boldsymbol{X}) = \boldsymbol{x}(\boldsymbol{X}) - \boldsymbol{X}
$$
The local deformation at a point is characterized by the {deformation gradient tensor} $\boldsymbol{F}:\Omega \longrightarrow \mathbb{R}^{2\times2}$, which is defined as the gradient of the spatial position with respect to the reference position:
$$
\boldsymbol{F} = \frac{\partial\boldsymbol{x}}{\partial\boldsymbol{X}} = \boldsymbol{I} + \nabla\boldsymbol{u}
$$
where $\boldsymbol{I}$ is the second-order identity tensor and $\nabla$ is the gradient operator with respect to $\boldsymbol{X}$.

From the deformation gradient, we can define objective (frame-indifferent) measures of stretching and rotation. The {right and left Cauchy-Green stretch tensors}, denoted by $\boldsymbol{C}$ and $\boldsymbol{B}$ respectively, are symmetric positive-definite tensors given by \cite{Gurtin1981}:
$$
\boldsymbol{C} = \boldsymbol{F}^{T}\boldsymbol{F} \quad \text{and} \quad \boldsymbol{B} = \boldsymbol{F}\boldsymbol{F}^{T}
$$
Based on these, we define two fundamental nonlinear strain measures. The {Green-St. Venant strain tensor} $\boldsymbol{E}:\Omega \longrightarrow Sym(\mathbb{R}^{2\times2})$ measures strain with respect to the reference configuration and is expressed as:
$$
\boldsymbol{E} = \frac{1}{2}(\boldsymbol{C}-\boldsymbol{I}) = \frac{1}{2}\left(\nabla\boldsymbol{u} + (\nabla\boldsymbol{u})^T + (\nabla\boldsymbol{u})^T\nabla\boldsymbol{u}\right)
$$
Conversely, the {Almansi-Hamel strain tensor} $\boldsymbol{e}$ measures strain with respect to the current configuration and is defined as:
$$
\boldsymbol{e} = \frac{1}{2}(\boldsymbol{I}-\boldsymbol{B}^{-1})
$$

In many engineering applications, the assumption of small deformations is valid, meaning the displacement gradients are very small ($\|\nabla\boldsymbol{u}\| \ll 1$). Under this assumption, the nonlinear terms in the Green-St. Venant tensor becomes negligible. Consequently, both $\boldsymbol{E}$ and $\boldsymbol{e}$ can be approximated by the same {linearized strain tensor} (or infinitesimal strain tensor) $\boldsymbol{\epsilon} \in Sym(\mathbb{R}^{2\times2})$, which is defined as the symmetric part of the displacement gradient:
$$
\boldsymbol{\epsilon} = \frac{1}{2}\left(\nabla\boldsymbol{u} + (\nabla\boldsymbol{u})^{T}\right)
$$
This simplification is central to the theories of linear elasticity and thermoelasticity.

\section{Rajagopal's implicit constitutive theory}

To overcome the limitations of classical linear elasticity, which predicts non-physical, unbounded strains at stress concentrations, our approach is founded on the implicit constitutive theory pioneered by Rajagopal \cite{rajagopal2007elasticity, rajagopal2011non}. This framework generalizes the classical stress-strain relationship by positing that the Cauchy stress tensor, $\boldsymbol{\sigma}$, and a suitable measure of deformation, such as the left Cauchy-Green stretch tensor, $\mathbf{B}$, are related implicitly through a constraint equation:
\begin{equation}
    \mathcal{F}(\mathbf{B}, \boldsymbol{\sigma}) = \boldsymbol{0}. \label{eq:rj1}
\end{equation}
This general form allows a much richer class of material responses than with traditional explicit models. In the above equation $\mathcal{F}(\cdot, \, \cdot)$ is the tensor-vlaued function. 

For the purposes of this study, which operates under the assumption of small displacement gradients, the general implicit relation in \eqref{eq:rj1} can be simplified. In this regime, it is convenient to work with an explicit relationship where the linearized strain, $\boldsymbol{\epsilon}$, is expressed as a nonlinear function of the Cauchy stress:
\begin{equation}
    \boldsymbol{\epsilon} = \mathcal{\widetilde{F}}(\boldsymbol{\sigma}). \label{eq:rj3}
\end{equation}
Here, $\mathcal{\widetilde{F}}: Sym(\mathbb{R}^{2\times2}) \rightarrow Sym(\mathbb{R}^{2\times2})$ is the nonlinear response function that uniquely defines the material's behavior.

A key feature of these models within this class is their ability to enforce a physical limit on deformation. A model is referred to as {strain-limiting} if the magnitude of the strain remains bounded for any admissible stress state \cite{MalliPhD2015, Mallikarjunaiah2015}. Mathematically, this condition is met if there exists a positive constant $M$ such that 
\[
\max_{\boldsymbol{\sigma} \in Sym(\mathbb{R}^{2\times2})} \|\mathcal{\widetilde{F}}(\boldsymbol{\sigma})\| \leq M.
\] 
Following recent successful applications in related areas \cite{Lee2022, yoon2021quasi, yoon2022MMS}, we adopt a specific form for this response function:
\begin{equation}
    \mathcal{\widetilde{F}}(\boldsymbol{\sigma}) = \frac{\mathbb{K}[\boldsymbol{\sigma}]}{\left(1 + (b\|\mathbb{K}^{1/2}[\boldsymbol{\sigma}]\|)^a\right)^{1/a}}. \label{eq:rj4}
\end{equation}
In this formulation, the numerator, $\mathbb{K}[\boldsymbol{\sigma}]$, represents the classical linear response governed by the fourth-order compliance tensor $\mathbb{K}$. The denominator acts as a nonlinear attenuation factor that becomes significant at high stress levels, ensuring the strain remains bounded. The positive constants $a$ and $b$ are modeling parameters: $a$ controls the sharpness of the transition from the linear to the limiting regime, while $b$ determines the ultimate strain limit, as 
\[
\sup_{\boldsymbol{\sigma} \in Sym(\mathbb{R}^{2\times2})} \|\mathcal{\widetilde{F}}(\boldsymbol{\sigma})\| \leq 1/b. 
\]
Notably, this nonlinear model gracefully reduces to the classical linear model ($\boldsymbol{\epsilon} = \mathbb{K}[\boldsymbol{\sigma}]$) in the limit as $b \to 0$ or $a \to \infty$, demonstrating that it is a direct generalization of the classical linear theory.

The compliance tensor $\mathbb{K}$ is the inverse of the fourth-order elasticity (or stiffness) tensor, $\mathbb{E}$. For the transversely isotropic materials considered in this work, the elasticity tensor, which maps strain to stress, is given by:
\begin{equation}
    \mathbb{E}[\boldsymbol{\epsilon}] = 2\mu\boldsymbol{\epsilon} + \lambda(\mathrm{tr}(\boldsymbol{\epsilon}))\mathbf{I} + \gamma(\boldsymbol{\epsilon}:\mathbf{M})\mathbf{M}. \label{eq:ts1}
\end{equation}
Here, $\lambda > 0$ and $\mu > 0$ are the classical Lamé parameters that describe the isotropic part of the material response. The material's anisotropy is captured by the additional parameter $\gamma$ and the {structural tensor} $\mathbf{M} = \mathbf{m} \otimes \mathbf{m}$, where $\mathbf{m}$ is a unit vector representing the preferred material direction (e.g., fiber orientation) \cite{Mallikarjunaiah2015,bhatta2025computational,ghosh2025computational,ghosh2025finite}.

For the proposed constitutive model to be physically meaningful and mathematically robust, its response function $\mathcal{\widetilde{F}}$ must satisfy several key properties. As established in \cite{itou2018states}, the chosen function meets the necessary criteria to guarantee that the resulting boundary value problem is well-posed, meaning a unique solution exists that depends continuously on the input data. These essential properties are outlined below.

The response function is characterized by four crucial mathematical conditions. First, it is {uniformly bounded}:
\begin{itemize}
    \item[(a)] \textbf{Boundedness:} $\|\mathcal{\widetilde{F}}(\mathbf{A})\| \leq \frac{1}{b}$ for all $\mathbf{A} \in Sym(\mathbb{R}^{2\times2})$. This property is the mathematical expression of the core ``strain-limiting'' feature, ensuring that the strain magnitude can never exceed a physical ceiling of $1/b$, regardless of the stress level.
\end{itemize}
Second, the function exhibits strict monotonicity, which ensures a direct and unambiguous relationship between stress and strain:
\begin{itemize}
    \item[(b)] \textbf{Strict Monotonicity:} $(\mathcal{\widetilde{F}}(\mathbf{A}_1) - \mathcal{\widetilde{F}}(\mathbf{A}_2)) : (\mathbf{A}_1 - \mathbf{A}_2) > 0$ for all distinct $\mathbf{A}_1, \mathbf{A}_2 \in Sym(\mathbb{R}^{2\times2})$. This implies that an increase in stress is always required to produce a further increase in strain. Mathematically, this property is critical as it guarantees the {uniqueness} of the solution to the boundary value problem.
\end{itemize}
Third and fourth, the function's continuity and growth are controlled by Lipschitz continuity and coercivity, which are vital for proving the existence of a solution and for ensuring the stability of numerical methods:
\begin{itemize}
    \item[(c)] \textbf{Lipschitz Continuity:} $\|\mathcal{\widetilde{F}}(\mathbf{A}_1) - \mathcal{\widetilde{F}}(\mathbf{A}_2)\| \leq \tilde{c}_1 \|\mathbf{A}_1 - \mathbf{A}_2\|$ for some constant $\tilde{c}_1 > 0$. This ensures that the strain response is continuous and does not change erratically with small changes in stress, a condition essential for the convergence of iterative numerical solvers.
    \item[(d)] \textbf{Coercivity:} $\langle \mathcal{\widetilde{F}}(\mathbf{A}), \mathbf{A} \rangle \geq \tilde{c}_2 \|\mathbf{A}\|^2$ for some constant $\tilde{c}_2 > 0$. This condition is related to the growth of the strain energy and is a technical requirement used in functional analysis to prove the existence of a solution within a variational framework.
\end{itemize}

Beyond these foundational properties, the model is also {hyperelastic}, as the constitutive law can be derived from a scalar energy potential. This is physically significant because it ensures that the material is path-independent and non-dissipative (i.e., no energy is lost during a closed loading-unloading cycle). Furthermore, for sufficiently large values of the nonlinearity parameter $b$ (equivalently, sufficiently small values of strains), the constitutive relation in \eqref{eq:rj3} is {invertible} \cite{mai2015monotonicity, mai2015strong}. This invertibility is crucial for implementation within a standard displacement-based finite element method. In such a framework, strains are first computed from the displacement field, and the corresponding stress tensor must then be determined by solving the inverse relation $\boldsymbol{\sigma} = \mathcal{\widetilde{F}}^{-1}(\boldsymbol{\epsilon})$ at each numerical integration point. While this step requires solving a local nonlinear problem, the well-defined structure of the model makes this a computationally tractable procedure.

\section{Constitutive model for strain-limiting thermoelasticity}

To model the coupled thermo-mechanical behavior, we extend the isothermal strain-limiting framework to include thermal effects. Following the general approach of Bustamante and Rajagopal for thermoelasticity \cite{bustamante2017}, the constitutive law in a fully nonlinear setting relates the Second Piola-Kirchhoff stress tensor $\boldsymbol{S}$, the Green-St. Venant strain tensor $\boldsymbol{E}$, and the temperature field $\theta$ through a general response function:
\begin{equation}
    \boldsymbol{S} = \mathcal{F}(\boldsymbol{E}, \theta). \label{eq:rj5}
\end{equation}
For the present study, which operates under the assumption of small displacement gradients, this relationship can be simplified. We adopt an additive decomposition where the total stress, denoted here as $\boldsymbol{\sigma}_{\text{Th}}$, is composed of a nonlinear mechanical stress component, $\boldsymbol{\sigma}$, and a thermal stress component arising from thermal expansion. In the context of small strains, this is expressed as:
\begin{equation}
    \boldsymbol{\sigma}_{\text{Th}} = \boldsymbol{\sigma} - \alpha\theta\mathbf{I}. \label{eq:ts2_revised}
\end{equation}
Here, $\alpha = \alpha_{T}(3\lambda + 2\mu)$, where $\alpha_T$ is the linear thermal expansion coefficient, and the term $(3\lambda + 2\mu)$ represents the bulk modulus for the isotropic part of the material.

The mechanical stress component, $\boldsymbol{\sigma}$, is governed by the strain-limiting constitutive law. While the model was initially formulated to express strain as a function of stress, its invertibility for sufficiently small strains \cite{mai2015monotonicity, mai2015strong, Mallikarjunaiah2015} allows us to express stress as a function of strain. This inverted form is essential for a displacement-based finite element formulation and is given by:
\begin{equation}
    \boldsymbol{\sigma}(\boldsymbol{\epsilon}) = \frac{\mathbb{E}[\boldsymbol{\epsilon}]}{\left(1 - (b\|\mathbb{E}^{1/2}[\boldsymbol{\epsilon}]\|)^a\right)^{1/a}}. \label{eq:sl2}
\end{equation}
 This formulation ensures that the mechanical stress response is inherently nonlinear and respects the strain-limiting principle, while the thermal effects are incorporated in a classically linear manner.

\subsection{Governing Equations and Problem Formulation}

In the absence of external body forces, the behavior of the coupled thermo-mechanical system is described by two governing partial differential equations. The first is the steady-state heat conduction equation for the temperature field $\theta$, and the second is the quasi-static momentum balance equation for the mechanical fields. The coupled system is given by:
\begin{align}
    -\nabla \cdot (k\nabla\theta) &= Q \quad \text{in} \quad \Omega, \label{eq:th1} \\
    -\nabla \cdot \boldsymbol{\sigma}_{\text{Th}} &= \boldsymbol{0} \quad \text{in} \quad \Omega. \label{eq:th2_revised}
\end{align}
Here, $k$ is the thermal conductivity of the material (assumed constant), and $Q:\Omega \to \mathbb{R}$ is an internal heat source or sink. Substituting the constitutive relation \eqref{eq:ts2_revised} into the momentum balance equation \eqref{eq:th2_revised} yields:
\begin{equation}
    -\nabla \cdot \boldsymbol{\sigma} + \alpha\nabla\theta = \boldsymbol{0} \quad \text{in} \quad \Omega. \label{eq:th2}
\end{equation}
This system exhibits a {one-way coupling}: the thermal field, through its gradient $\nabla\theta$, acts as a body force on the mechanical system, but the mechanical deformation does not influence the temperature distribution. This decoupling allows for a convenient and efficient sequential solution strategy.

\subsection{Sequential Solution Strategy}

Based on the one-way coupling, we first solve for the temperature field and then use that to solve for the mechanical displacement.

\paragraph{Step 1: Solve the Thermal Boundary Value Problem.}
The temperature distribution $\theta$ across the domain $\Omega$ is obtained by solving the linear, second-order elliptic PDE:
\begin{equation}
    -\nabla \cdot (k\nabla\theta) = Q \quad \text{in} \quad \Omega, \label{eq:th3}
\end{equation}
subject to a combination of Dirichlet and Neumann boundary conditions on complementary parts of the boundary $\partial\Omega = \overline{\partial\Omega_{D}^{\theta}} \cup \overline{\partial\Omega_{N}^{\theta}}$:
\begin{align*}
    \theta &= \theta_I \quad \text{on} \quad \partial\Omega_{D}^{\theta} \quad \text{(Prescribed Temperature)}, \\
    \boldsymbol{n} \cdot (k\nabla\theta) &= q_0 \quad \text{on} \quad \partial\Omega_{N}^{\theta} \quad \text{(Prescribed Heat Flux)}.
\end{align*}

\paragraph{Step 2: Solve the Mechanical Boundary Value Problem.}
Once the temperature field $\theta$ is known, its gradient is computed and treated as a known thermal body force, $\mathbf{f} = -\alpha\nabla\theta$. The displacement field $\mathbf{u}$ is then found by solving the following quasi-linear elliptic system, which is obtained by substituting the constitutive law \eqref{eq:sl2} into the momentum balance equation \eqref{eq:th2}:
\begin{equation}
    -\nabla \cdot \left(\frac{\mathbb{E}[\boldsymbol{\epsilon}(\mathbf{u})]}{\left(1 - (b\|\mathbb{E}^{1/2}[\boldsymbol{\epsilon}(\mathbf{u})]\|)^a\right)^{1/a}}\right) = \mathbf{f} \quad \text{in} \quad \Omega, \label{eq:th4}
\end{equation}
subject to appropriate mechanical boundary conditions on $\partial\Omega = \overline{\partial\Omega_D} \cup \overline{\partial\Omega_N}$:
\begin{align*}
    \mathbf{u} &= \mathbf{u}_I \quad \text{on} \quad \partial\Omega_D \quad \text{(Prescribed Displacement)}, \\
    \boldsymbol{\sigma}\mathbf{n} &= \mathbf{t}_0 \quad \text{on} \quad \partial\Omega_N \quad \text{(Prescribed Traction)}.
\end{align*}

\subsection{Well-posedness of the Mechanical Governing Equations}

This section is dedicated to establishing the mathematical well-posedness of the nonlinear mechanical boundary value problem. A problem is considered well-posed if a solution exists, is unique, and depends continuously on the input data. We demonstrate the existence and uniqueness of a weak solution by constructing a suitable variational (or weak) formulation and leveraging key properties of the strain-limiting constitutive model. The theoretical approach closely follows the functional analysis framework introduced by Beck et al. \cite{Beck2017} for implicit constitutive theories.

\subsubsection{Strong and Weak Formulations}

The mechanical response of the body is governed by the quasi-static equilibrium equation coupled with the nonlinear strain-limiting constitutive law and appropriate boundary conditions. This constitutes the strong form of the boundary value problem (BVP). For our specific case, the body force $\mathbf{f}$ is determined by the thermal field ($\mathbf{f} = -\alpha\nabla\theta$), and we consider a traction-free Neumann boundary ($\mathbf{g} = \mathbf{0}$). The BVP is thus stated as finding the displacement $\mathbf{u}$ and stress $\boldsymbol{\sigma}$ such that:
\begin{equation}
\begin{cases}
    -\nabla \cdot \boldsymbol{\sigma} = \mathbf{f} & \text{in } \Omega, \\
    \boldsymbol{\epsilon}(\mathbf{u}) = \mathcal{\widetilde{F}}(\boldsymbol{\sigma}) := \frac{\mathbb{K}[\boldsymbol{\sigma}]}{(1+(b\|\mathbb{K}^{1/2}[\boldsymbol{\sigma}]\|)^{a})^{1/a}} & \text{in } \Omega, \\
    \mathbf{u} = \mathbf{u}_I & \text{on } \partial\Omega_D, \\
    \boldsymbol{\sigma}\mathbf{n} = \mathbf{g} & \text{on } \partial\Omega_N.
\end{cases}
\label{eq:bvp1}
\end{equation}
To analyze this nonlinear system, we reformulate it in a weak (or variational) form. This is achieved by multiplying the equilibrium equation by a suitable test function $\mathbf{w}$, integrating over the domain $\Omega$, and applying the divergence theorem. This process transfers a derivative from the stress tensor to the test function and naturally incorporates the Neumann boundary condition. The resulting weak formulation is to find a displacement field $\mathbf{u}$ that satisfies the boundary conditions and the following integral equation for all admissible test functions $\mathbf{w}$:
\begin{equation}
    \int_{\Omega} \boldsymbol{\sigma}(\boldsymbol{\epsilon}(\mathbf{u})) : \boldsymbol{\epsilon}(\mathbf{w}) \,d\mathbf{x} = \int_{\Omega} \mathbf{f} \cdot \mathbf{w} \,d\mathbf{x} + \int_{\partial\Omega_N} \mathbf{g} \cdot \mathbf{w} \,d\mathbf{S}. \label{eq:bvp3}
\end{equation}
The test functions $\mathbf{w}$ belong to a space of kinematically admissible virtual displacements, which vanish on the Dirichlet boundary $\partial\Omega_D$.

\subsubsection{Existence and Uniqueness of the Solution}

The proof of existence and uniqueness for the weak formulation \eqref{eq:bvp3} relies on the following set of physically and mathematically motivated assumptions:
\begin{enumerate}
    \item \textbf{Material Model Regularity:} The model parameters $a, b, \lambda,$ and $\mu$ are assumed to be positive constants. This ensures the elasticity tensor is positive definite. Furthermore, it is required that in the limit of vanishing nonlinearity ($b \to 0^+$), the strain-limiting model consistently reduces to the classical linear elastic model.
    \item \textbf{Global Equilibrium:} In the specific case of a pure traction problem (where the Dirichlet boundary $\partial\Omega_D$ is empty), the externally applied forces must be in static equilibrium to preclude rigid body motion. This requires that the net force on the body is zero:
    \begin{equation}
        \int_{\Omega} \mathbf{f} \,d\mathbf{x} + \int_{\partial\Omega_N} \mathbf{g} \,d\mathbf{S} = \mathbf{0}. \label{eq:bvp2}
    \end{equation}
    \item \textbf{Data Regularity:} The prescribed boundary displacement $\mathbf{u}_I$ is assumed to possess sufficient regularity, specifically $\mathbf{u}_I \in (W^{1,1}(\Omega))^2$. Additionally, the strain resulting from this prescribed displacement, $\boldsymbol{\epsilon}(\mathbf{u}_I(\mathbf{x}))$, must be contained within a compact set in $\mathbb{R}^{2\times2}$ for almost every $\mathbf{x} \in \overline{\Omega}$. These are technical conditions on the smoothness of the boundary data required by the proof.
\end{enumerate}

With these assumptions in place, we can state the main theorem regarding the well-posedness of the problem.

\vspace{1em}
\noindent
\textbf{Theorem.} \textit{Let $\Omega \subset \mathbb{R}^2$ be an open, bounded, and connected Lipschitz domain with a boundary $\partial\Omega$ partitioned into Dirichlet ($\partial\Omega_D$) and Neumann ($\partial\Omega_N$) segments. Given a body force $\mathbf{f} \in (L^2(\Omega))^2$, a traction vector $\mathbf{g} \in (L^2(\partial\Omega_N))^2$, and a prescribed displacement $\mathbf{u}_I$ that satisfy the preceding assumptions, the boundary value problem defined by \eqref{eq:bvp1} is well-posed. Specifically, there exists a unique weak solution pair $(\mathbf{u}, \boldsymbol{\sigma})$ in the function space $(W^{1,1}(\Omega))^2 \times Sym(L^1(\Omega))^{2\times2}$ that satisfies the weak formulation \eqref{eq:bvp3}. Existence is guaranteed by the coercive and continuous nature of the underlying energy functional, while uniqueness is a direct consequence of the strict monotonicity of the constitutive operator $\mathcal{\widetilde{F}}$.}

\section{Numerical Solution using the Finite Element Method}

The coupled, nonlinear system of partial differential equations is solved using the Galerkin finite element method \cite{Ciarlet2002}. The solution process begins by recasting the strong form of the governing equations into a more general weak (or variational) formulation, which is then discretized over a computational mesh. The resulting nonlinear system of algebraic equations is subsequently solved using an iterative linearization scheme.

\subsection{Variational Formulation and Function Spaces}

The foundation of the finite element method is the variational formulation of the boundary value problem. To formally define this, we first introduce the necessary function spaces. Let $L^2(\Omega)$ denote the space of square-integrable functions on the domain $\Omega$, equipped with the standard inner product $(u,v) = \int_{\Omega} uv \,d\mathbf{x}$ and the associated norm $\|u\|_{L^2} = (u,u)^{1/2}$. For problems involving second-order differential operators, the appropriate mathematical setting is the Sobolev space $H^1(\Omega)$, defined as:
\[
H^1(\Omega) = \{ v \in L^2(\Omega) : \nabla v \in (L^2(\Omega))^2 \}
\]
Using this space, we define the trial and test spaces for the temperature and displacement fields. The trial (or solution) spaces contain functions that satisfy the essential (Dirichlet) boundary conditions, while the test (or variation) spaces contain functions satisfying the corresponding \textit{homogeneous} boundary conditions.

For the temperature field, the trial and test spaces are:
\begin{align*}
    V^{\theta} &= \{ \theta \in H^1(\Omega) : \theta = \theta_I \text{ on } \partial\Omega_D^{\theta} \}, \\
    V_0^{\theta} &= \{ q \in H^1(\Omega) : q = 0 \text{ on } \partial\Omega_D^{\theta} \}.
\end{align*}
Similarly, for the vector-valued displacement field, the spaces are:
\begin{align*}
    \mathbf{V} &= \{ \mathbf{u} \in (H^1(\Omega))^2 : \mathbf{u} = \mathbf{u}_I \text{ on } \partial\Omega_D \}, \\
    \mathbf{V}_0 &= \{ \mathbf{v} \in (H^1(\Omega))^2 : \mathbf{v} = \mathbf{0} \text{ on } \partial\Omega_D \}.
\end{align*}

The continuous weak formulation is derived by multiplying the strong-form equations \eqref{eq:th3} and \eqref{eq:th4} by arbitrary test functions ($q \in V_0^{\theta}$ and $\mathbf{v} \in \mathbf{V}_0$, respectively) and integrating over the domain $\Omega$. Applying the divergence theorem (integration by parts) yields the following variational problem:

Find $(\theta, \mathbf{u}) \in V^\theta \times \mathbf{V}$ such that:
\begin{align}
    a_{\theta}(\theta, q) &= l_{\theta}(q) \quad \forall q \in V_0^{\theta}, \label{eq:wk1} \\
    a(\mathbf{u}; \mathbf{v}) &= l(\theta; \mathbf{v}) \quad \forall \mathbf{v} \in \mathbf{V}_0. \label{eq:wk2}
\end{align}
The thermal problem is linear, with the bilinear form $a_{\theta}(\cdot, \cdot)$ and linear functional $l_{\theta}(\cdot)$ defined as:
\begin{align}
    a_{\theta}(\theta, q) &= \int_{\Omega} k \nabla\theta \cdot \nabla q \, d\Omega, \label{eq:wk3} \\
    l_{\theta}(q) &= \int_{\Omega} Qq \, d\Omega. \label{eq:wk4}
\end{align}
The mechanical problem is nonlinear. We define its semi-linear form $a(\cdot; \cdot)$ and linear functional $l(\cdot; \cdot)$ as:
\begin{align}
    a(\mathbf{u}; \mathbf{v}) &= \int_{\Omega} \frac{\mathbb{E}[\boldsymbol{\epsilon}(\mathbf{u})] : \boldsymbol{\epsilon}(\mathbf{v})}{\left(1 - (b\|\mathbb{E}^{1/2}[\boldsymbol{\epsilon}(\mathbf{u})]\|)^a\right)^{1/a}} \, d\Omega, \label{eq:wk5} \\
    l(\theta; \mathbf{v}) &= -\int_{\Omega} (\alpha\nabla\theta) \cdot \mathbf{v} \, d\Omega. \label{eq:wk6}
\end{align}
Note that the form $a(\mathbf{u}; \mathbf{v})$ is nonlinear with respect to its first argument, $\mathbf{u}$, which necessitates an iterative solution approach.

\subsection{Finite Element Discretization}

The continuous variational problem posed in infinite-dimensional spaces is made computationally tractable by discretization. The domain $\overline{\Omega}$ is partitioned into a mesh, $\mathcal{T}_h$, consisting of a set of non-overlapping quadrilateral elements $K$. The parameter $h = \max_{K \in \mathcal{T}_h} (\text{diam}(K))$ denotes the mesh size.

On this mesh, we construct finite-dimensional subspaces of the continuous function spaces. Let $Q_k(K)$ be the space of polynomials of degree up to $k$ in each variable on an element $K$. We define the discrete spaces for the temperature and displacement fields as:
\begin{align*}
    S_h^\theta &= \{ q_h \in C^0(\overline{\Omega}) : q_h|_K \in Q_k(K) \quad \forall K \in \mathcal{T}_h \}, \\
    \mathbf{S}_h &= \{ \mathbf{u}_h \in (C^0(\overline{\Omega}))^2 : \mathbf{u}_h|_K \in (Q_k(K))^2 \quad \forall K \in \mathcal{T}_h \}.
\end{align*}
The discrete trial and test spaces are then defined by intersecting these polynomial spaces with the appropriate continuous function spaces:
\begin{align*}
    V_h^\theta = S_h^\theta \cap V^\theta, \quad V_{h,0}^\theta = S_h^\theta \cap V_0^\theta, \\
    \mathbf{V}_h = \mathbf{S}_h \cap \mathbf{V}, \quad \mathbf{V}_{h,0} = \mathbf{S}_h \cap \mathbf{V}_0.
\end{align*}
The discrete weak formulation is obtained by restricting the continuous problem to these finite-dimensional subspaces.

\subsection{Solution of the Nonlinear System via Picard Iteration}

Due to the one-way coupling, the problem is solved sequentially. First, the linear discrete thermal problem is solved for the temperature field $\theta_h$. Then, the nonlinear discrete mechanical problem is solved for the displacement $\mathbf{u}_h$.

The nonlinearity in the mechanical weak form \eqref{eq:wk5} is handled using {Picard's iterative method}, a fixed-point linearization technique. The method generates a sequence of solutions $\{\mathbf{u}_h^n\}_{n=0}^\infty$ that converges to the true solution. The core idea is to linearize the problem at each step $n+1$ by evaluating the nonlinear denominator using the solution from the previous step, $\mathbf{u}_h^n$.

The iterative algorithm proceeds as follows:
\begin{enumerate}
    \item \textbf{Initialization (n=0):} Compute an initial guess, $\mathbf{u}_h^0$, by solving the corresponding linear elasticity problem, which is achieved by setting the nonlinearity parameter $b=0$.
    \item \textbf{Iteration (n+1):} For $n=0, 1, 2, \dots$, given the previous solution $\mathbf{u}_h^n$, find the next iterate $\mathbf{u}_h^{n+1} \in \mathbf{V}_h$ by solving the following \textit{linear} variational problem for all $\mathbf{v}_h \in \mathbf{V}_{h,0}$:
    \begin{equation}
        \int_{\Omega} \frac{\mathbb{E}[\boldsymbol{\epsilon}(\mathbf{u}_h^{n+1})] : \boldsymbol{\epsilon}(\mathbf{v}_h)}{\left(1 - (b\|\mathbb{E}^{1/2}[\boldsymbol{\epsilon}(\mathbf{u}_h^n)]\|)^a\right)^{1/a}} \, d\Omega = -\int_{\Omega} (\alpha\nabla\theta_h) \cdot \mathbf{v}_h \, d\Omega. \label{eq:wk7}
    \end{equation}
    \item \textbf{Convergence Check:} The iteration is terminated when the difference between successive solutions is sufficiently small in a suitable norm, e.g., when $\|\mathbf{u}_h^{n+1} - \mathbf{u}_h^n\|_{L^2} < \text{tol}$, where `tol` is a prescribed tolerance.
\end{enumerate}

\section{Numerical Results and Discussion}

This section presents a series of numerical simulations designed to demonstrate the performance and validate the predictions of the proposed strain-limiting model for transversely isotropic bodies under thermo-mechanical loading. We investigate a canonical fracture mechanics problem: a rectangular plate containing a sharp edge crack. This configuration is specifically chosen because it induces a stress singularity at the crack tip within the framework of classical linear elasticity, providing a challenging test case to highlight the key features of our nonlinear model. The primary objective is to demonstrate that the strain-limiting constitutive law effectively regularizes the strain field, preventing the non-physical, infinite strain predicted by the linear model and yielding a more realistic stress and strain distribution in the high-stress region near the crack tip.

The governing equations are discretized and solved using a conventional continuous Galerkin finite element method. The computational domain is meshed with bilinear quadrilateral elements ($Q_2$ elements). All numerical experiments are implemented using the advanced, open-source \textsf{C++} finite element library \textsf{deal.II} \cite{Arndt2023dealii, Arndt2021dealii}, which provides a robust and flexible platform for solving complex partial differential equations.

The material nonlinearity introduced by the strain-limiting constitutive model is resolved using the Picard iterative procedure described in the previous section. To rigorously control the iterative process, we monitor the convergence at each step. The iteration is deemed to have converged when the solution stabilizes, which is quantified by computing the discrete $L^2$ norm of the difference between the displacement solution vectors from two consecutive iterations, $\|\mathbf{u}_h^{n+1} - \mathbf{u}_h^n\|_{L^2}$. The process is terminated once this norm falls below a prescribed tolerance, ensuring a robust and accurate final solution.  The stress and strain fields are post-processed to analyze the model's behavior.

\subsection{Setup of domain and boundary conditions:}

For the numerical simulations, we consider a rectangular material
with a single edge crack subjected to thermal and mechanical loading.
Our primary objective is to compare the stress, strain, and strain
energy density predicted by the strain-limiting model with the predictions
obtained from the linearized elastic model. For the numerical simulations,
Picard's iteration was used with a tolerance of $TOL=10^{-8}$ and
maximum iterations, $N_{max}=100$. Figure 1 illustrates the geometric
setup and boundary conditions for this problem. Here we use a rectangular
plate as our computational domain embedded with a horizontal crack
extending from the left edge along (parallel to the x-axis)
a line $y=0.5$, $0.5\le x\le1$

\begin{figure}[H]
\begin{centering}
\includegraphics[scale=0.4]{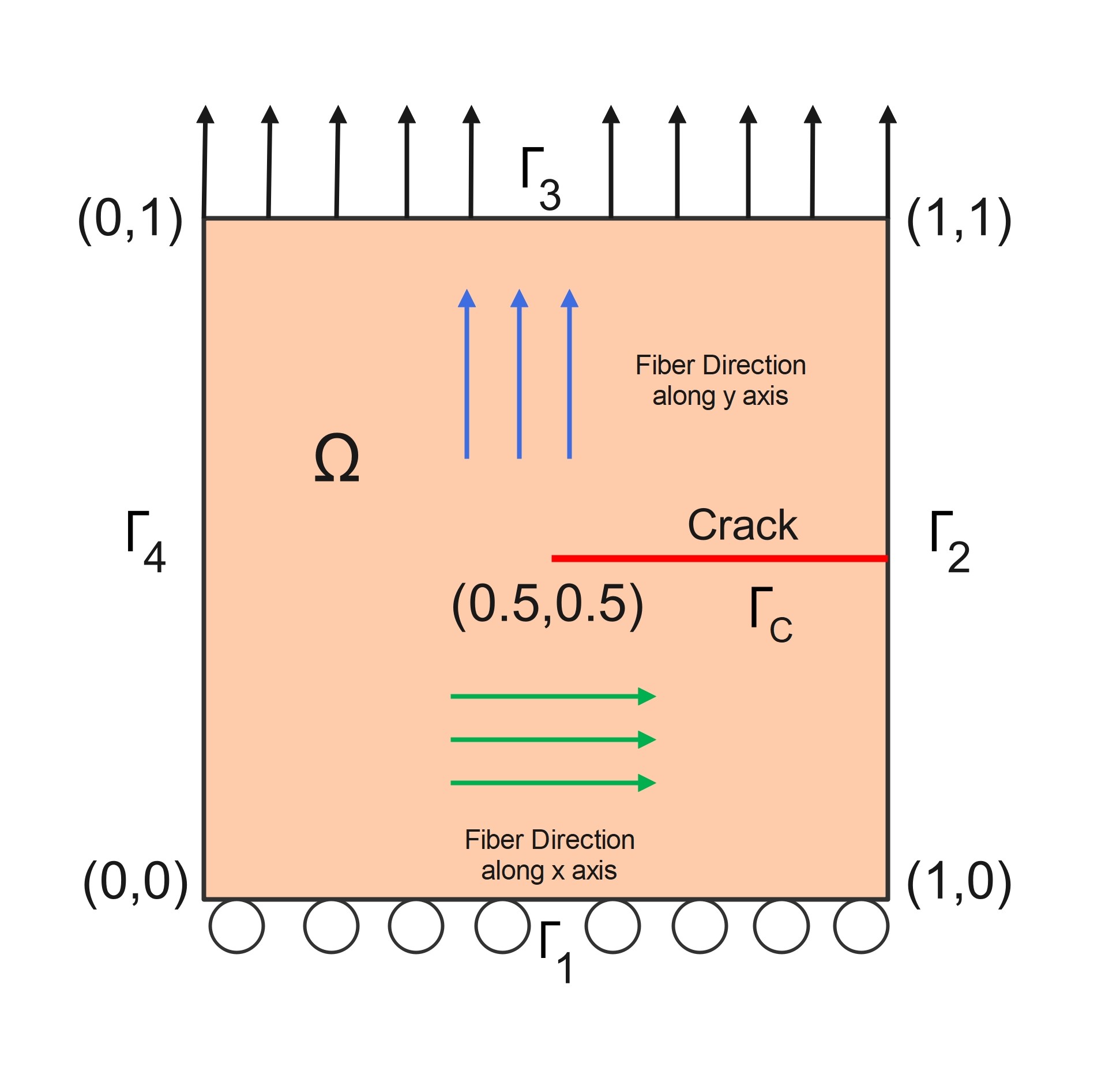}
\par\end{centering}
\centering{}\caption{Computational domain}
\end{figure}

Dirichlet boundary conditions are :
\begin{eqnarray*}
\boldsymbol{u}=(0,d) & \mathrm{on}\;\varGamma_{3}(\mathrm{top\;boundary)} & \mathrm{for\;displacement}\\
u_{y}=0 & \mathrm{on}\;\varGamma_{1}(\mathrm{bottom\;boundary)} & \mathrm{for\;displacement}\\
\theta=\theta_{0} & \mathrm{on}\;\varGamma_{1}(\mathrm{bottom\;boundary)} & \mathrm{for\;temperature}
\end{eqnarray*}

Neumann boundary conditions are :
\begin{eqnarray*}
\boldsymbol{\sigma n}=\boldsymbol{0} & \mathrm{on}\; \varGamma_{c},\varGamma_{2},\varGamma_{4}(\mathrm{right\;and\;left\;boundary)} & \mathrm{for\;displacement}\\
\nabla\theta\cdot\boldsymbol{n}=0 & \mathrm{on}\;\varGamma_{2},\varGamma_{3},\varGamma_{4}(\mathrm{right,top\;and\;left\;boundary)} & \mathrm{for\;temperature}
\end{eqnarray*}

\subsection{Fiber orientation parallel to the crack}

In this study, we analyze the scenario where the reinforcing fibers are oriented parallel to the crack, which lies along the x-axis. The material is assumed to exhibit {transverse isotropy}, meaning its mechanical and thermal properties are uniform in any direction on a given plane (the plane of isotropy) but differ along the axis perpendicular to that plane. This directional dependence is mathematically captured using the structural tensor $\mathbf{M}$. For fibers aligned with the x-axis, this tensor is defined as $\mathbf{M} = \boldsymbol{e}_{1} \otimes \boldsymbol{e}_{1}$, where $\boldsymbol{e}_{1}$ is the unit vector in the x-direction.

To evaluate the material's thermoelastic response, we investigate two distinct thermal loading scenarios applied to the bottom boundary of the computational domain:

\begin{enumerate}
    \item \textbf{Uniform Thermal Load}: A constant temperature is applied via a uniform Dirichlet boundary condition, $\theta = 100$. This case simulates a consistent, steady-state heating environment across the boundary.
    \item \textbf{Non-Uniform Thermal Gradient}: A spatially varying temperature profile is imposed through a parabolic Dirichlet boundary condition, given by $\theta(x) = 400x(1-x)$. This setup introduces a significant thermal gradient, mimicking more complex real-world conditions where heat distribution is uneven.
\end{enumerate}

By comparing the results from these two cases, we can conduct a comprehensive assessment of the thermoelastic behavior near the crack tip under both uniform and gradient-driven thermal loads.

\subsubsection{Constant temperature}

\begin{figure}[H]
\begin{centering}
\includegraphics[scale=0.2]{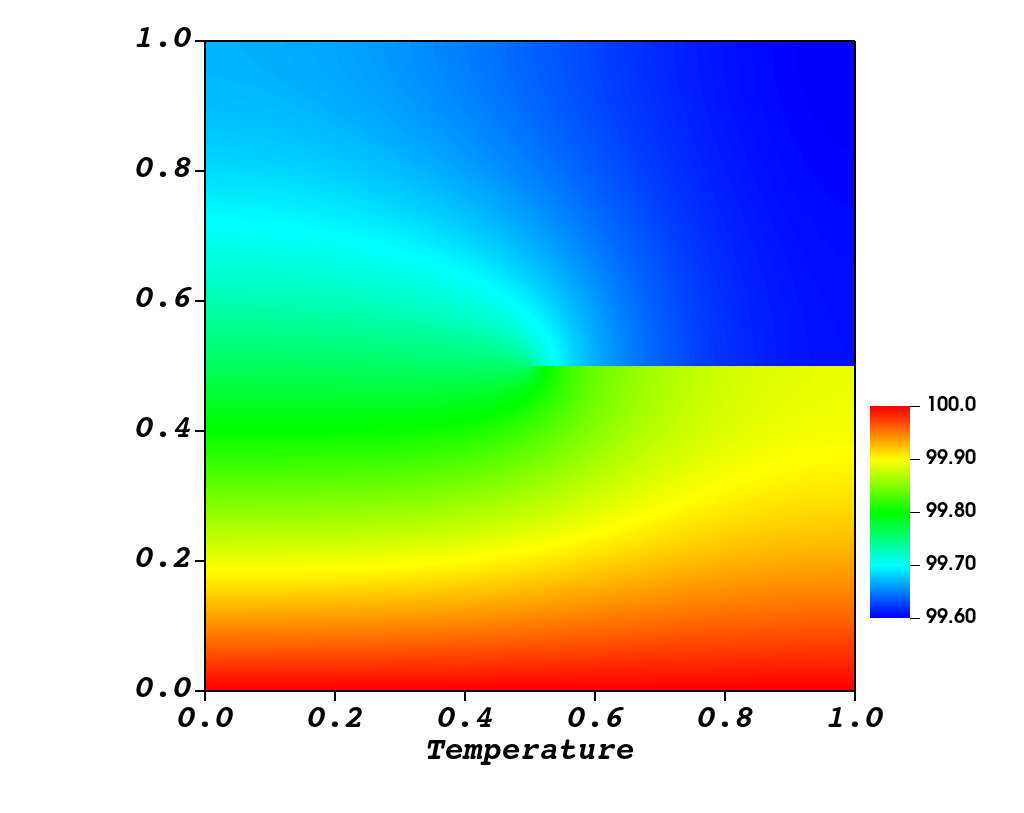}
\par\end{centering}
\centering{}\caption{Temperature Profile, $\theta=100,a=0.5,b=0.02$}\label{fig:Temp-const}
\end{figure}

In this example, a constant temperature $(\theta=100)$ boundary condition
is provided at the bottom boundary. We observe a very slight change
in the temperature (lowest = 99.60, highest =100) throughout the domain.

\begin{figure}[H]
\centering 
\includegraphics[scale=0.3]{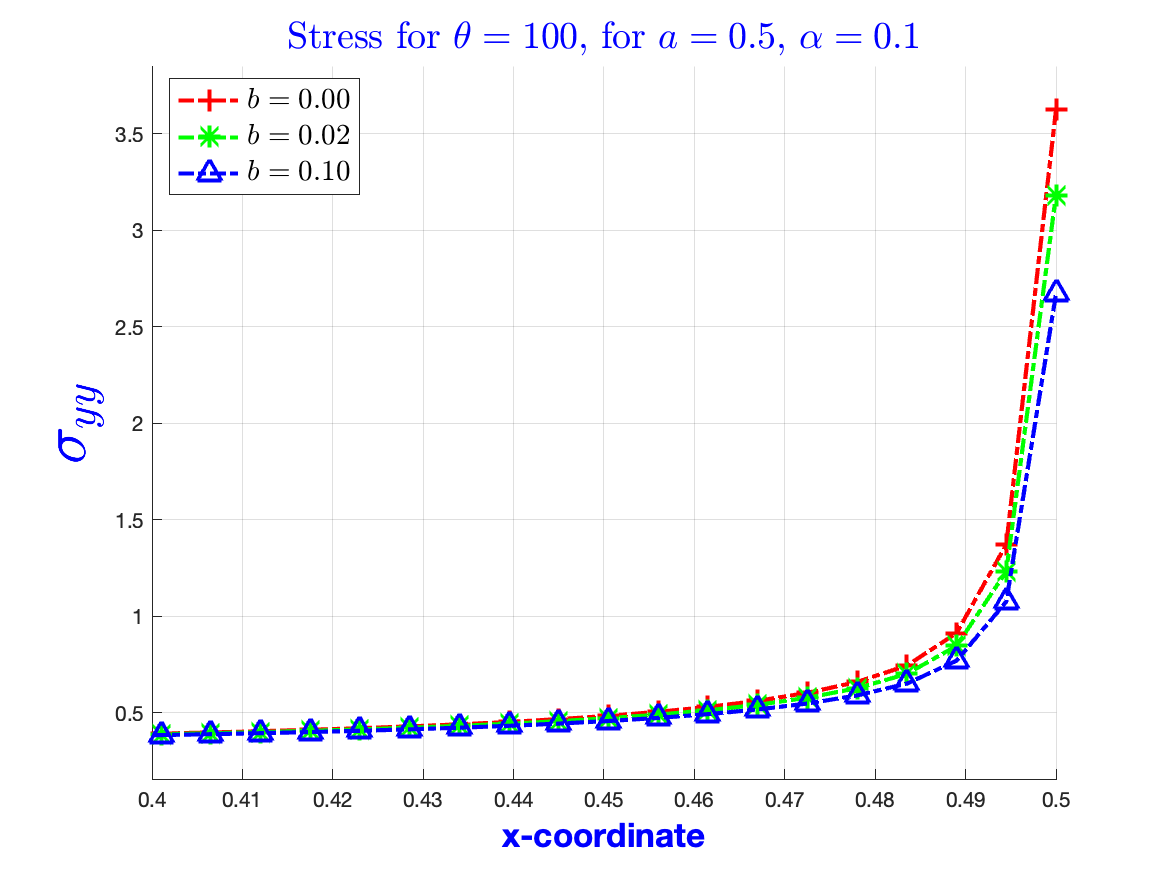}\quad
\includegraphics[scale=0.3]{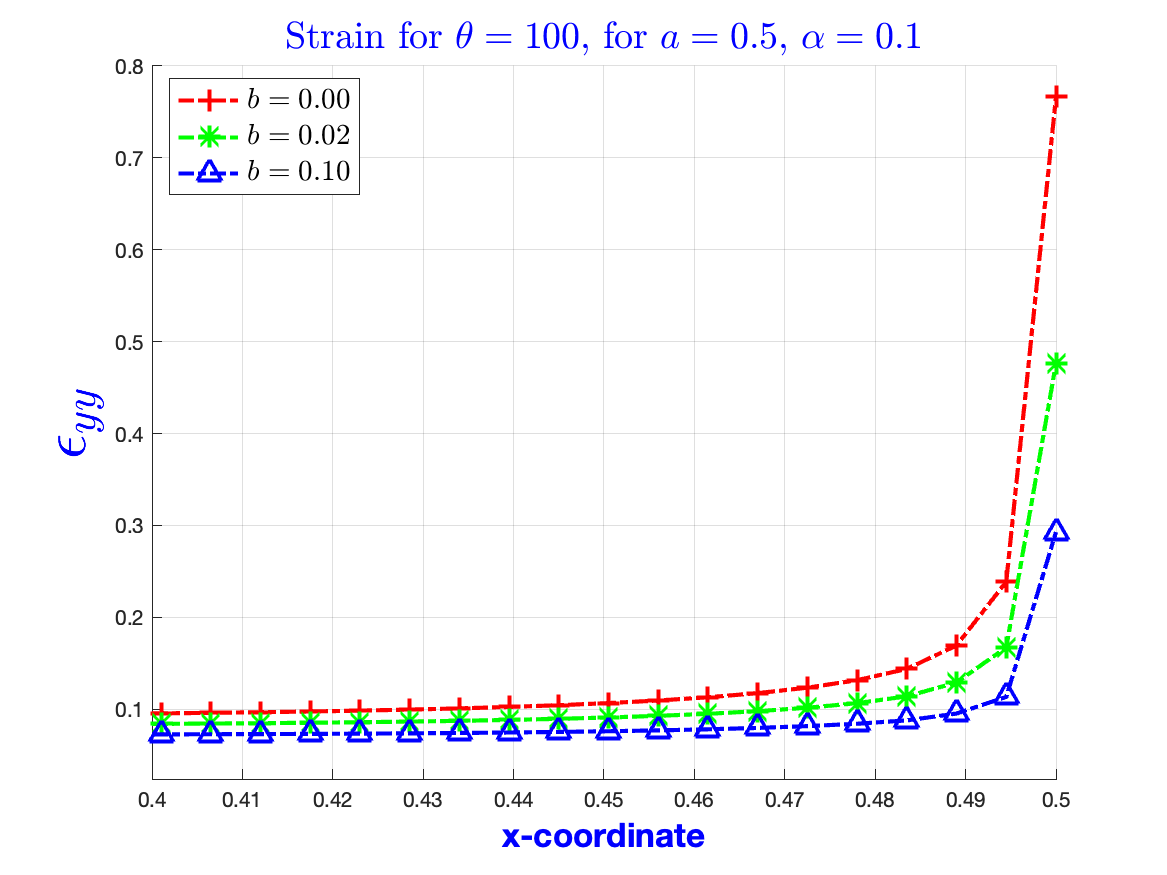}
\caption{Stress and strain for $b$ variation, $\theta=100,a=0.5$}\label{fig:SS_m1_bv_const}
\end{figure}

Figure (\ref{fig:SS_m1_bv_const}) illustrates the effect of the material parameter $b$ on the thermo-mechanical response of the material. A clear trend emerges from the data: as the value of $b$ increases, there is a corresponding moderate decrease in the peak stress and strain values observed near the crack tip. This inverse relationship suggests that the parameter $b$, which governs the material's non-linear behavior, plays a crucial role in mitigating deformation. In essence, higher values of $b$ enhance the material's ability to resist deformation under thermal load, thereby reducing stress concentration at the crack.

The maximum stress and strain values, recorded at {3.91996} and {0.794503}, respectively, occur specifically when $b=0$. This scenario corresponds to the classical linear thermoelastic model, where the non-linear mitigating effects are absent. Therefore, the linear case represents the baseline or "worst-case" scenario in terms of mechanical response at the crack tip.

\begin{figure}[H]
\centering 
\includegraphics[scale=0.3]{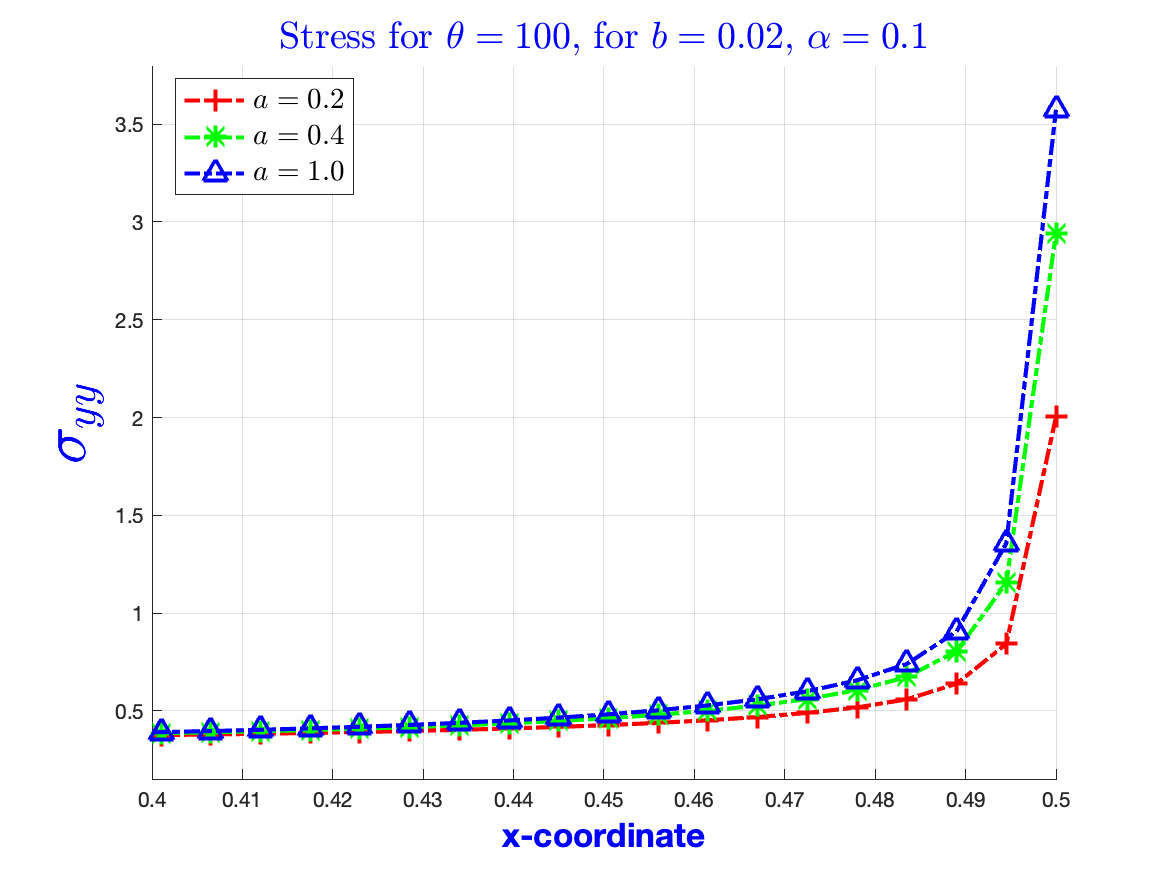}\quad
\includegraphics[scale=0.3]{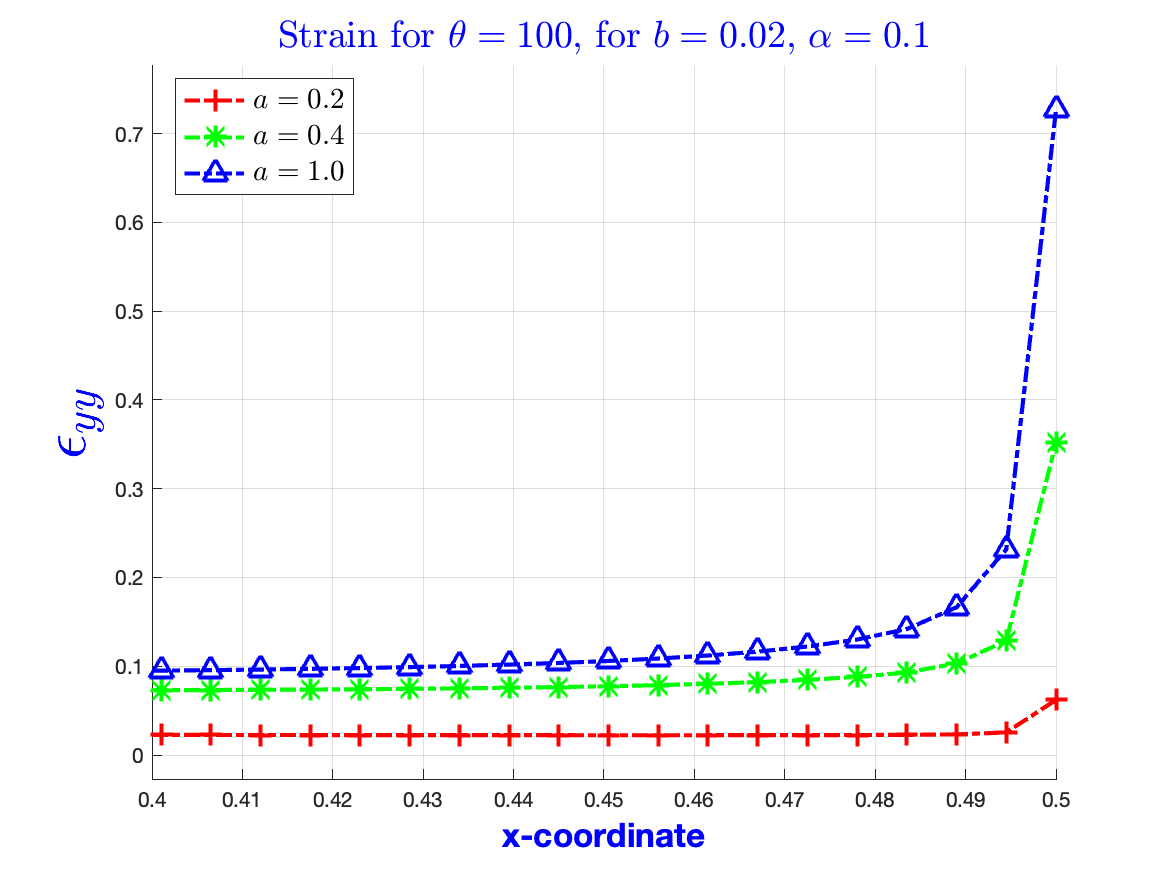}
\caption{Stress and strain for $a$ variation , $\theta=100,b=0.02$}\label{fig:SS_m1_av_const}
\end{figure}

Analysis of Figure (\ref{fig:SS_m1_av_const}) reveals a strong positive correlation between the material parameter $a$ and the system's mechanical response. Specifically, an increase in the value of $a$ leads to a significant amplification of both stress and strain within the material. This trend culminates when $a=1.0$, the upper limit considered in this analysis, at which point the maximum stress of {3.86449} and maximum strain of {0.752424} are observed.

This direct relationship has a critical implication for material performance: higher values of the parameter $a$ correspond to a lower resilience against external forces, indicating a heightened susceptibility to damage and potential fracture. Therefore, from a design and safety perspective, the parameter $a$ must be carefully controlled and ideally minimized to ensure the structural integrity and durability of the material under operational loads.

\begin{figure}[H]
\centering 
\includegraphics[scale=0.2]{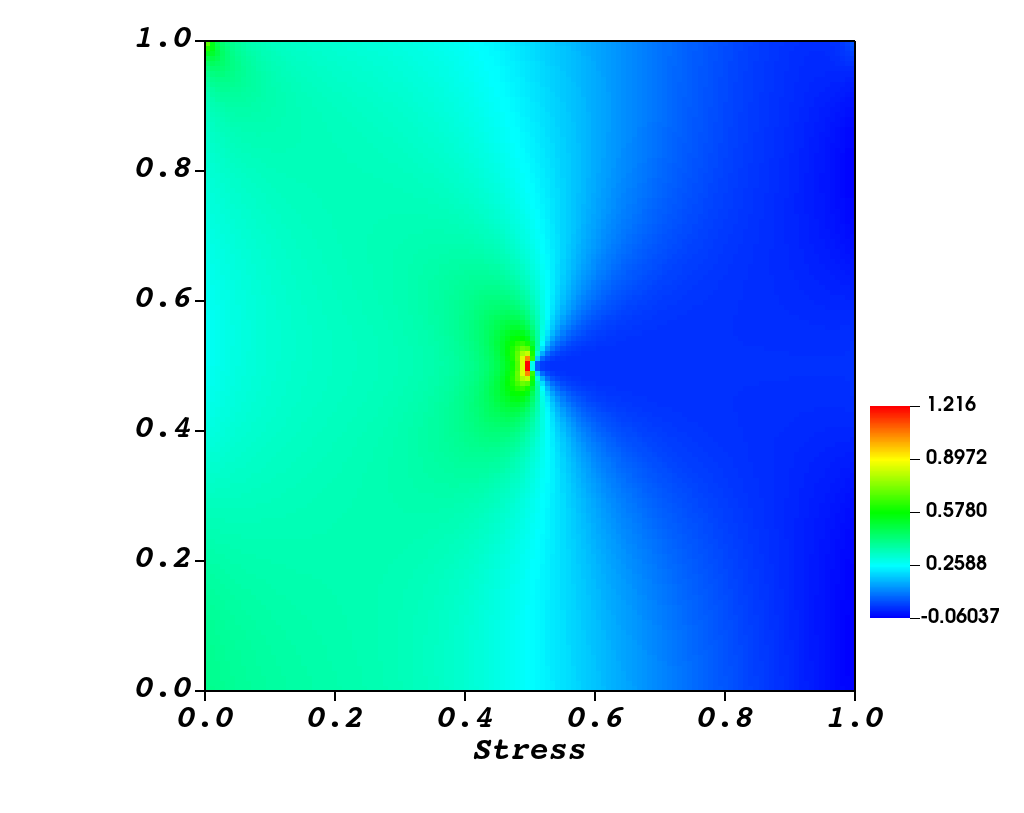}\quad
\includegraphics[scale=0.2]{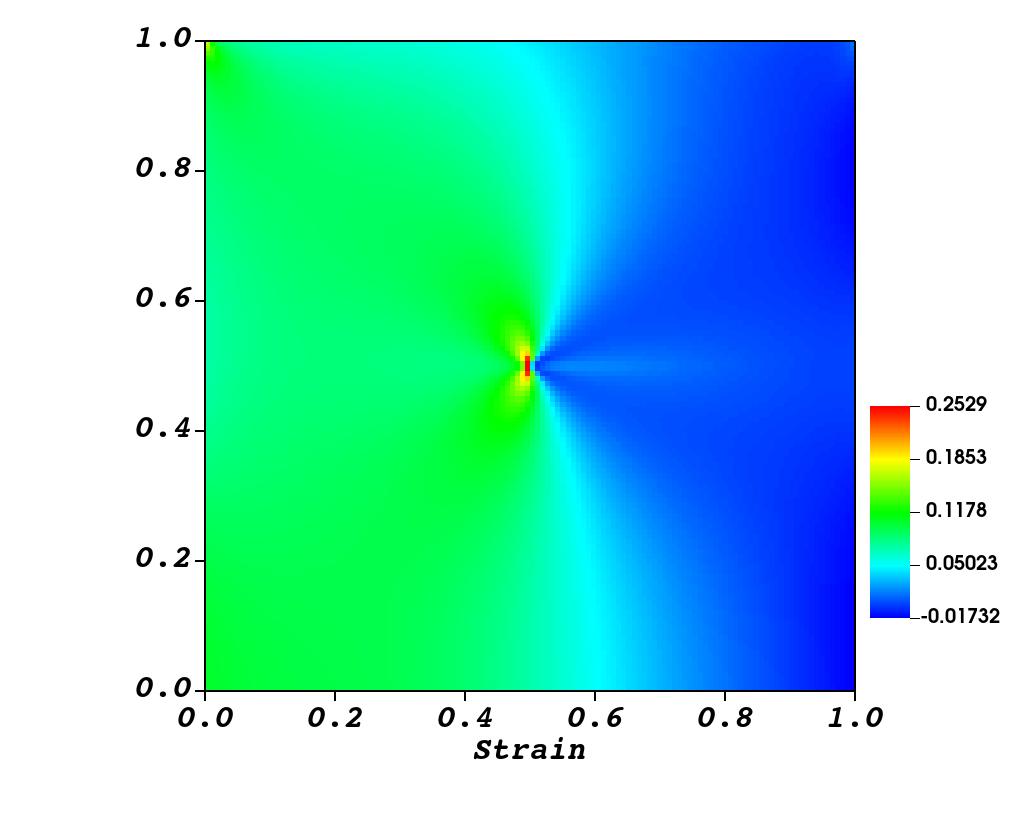}
\caption{Stress and strain for $\theta=100,a=0.5,b=0.02$}\label{fig:SS_m1_const_visit}
\end{figure}

Figure (\ref{fig:SS_m1_const_visit}) presents the distribution of axial stress and strain throughout the domain under a uniform thermal load, where a constant temperature of $\theta=100$ is applied to the bottom boundary. As is characteristic of fracture mechanics problems, a significant {stress and strain concentration} is observed in the immediate vicinity of the crack tip. The contour plot's sharp color gradient clearly delineates this high-stress zone from the rest of the body. Moving away from the crack, particularly into the right half of the domain, the stress and strain values diminish substantially.

Quantitatively, the analysis reveals that the peak tensile stress reaches a value of {1.216} at the crack tip, with a corresponding maximum strain of {0.2529}. The simulation also captures areas of slight compression elsewhere in the body, indicated by a minimum stress value of {-0.06037} and a minimum strain of {-0.01732}. This distribution underscores the crack's role as the primary driver of mechanical response under these thermal conditions.

\begin{figure}[H]
\centering 
\includegraphics[scale=0.2]{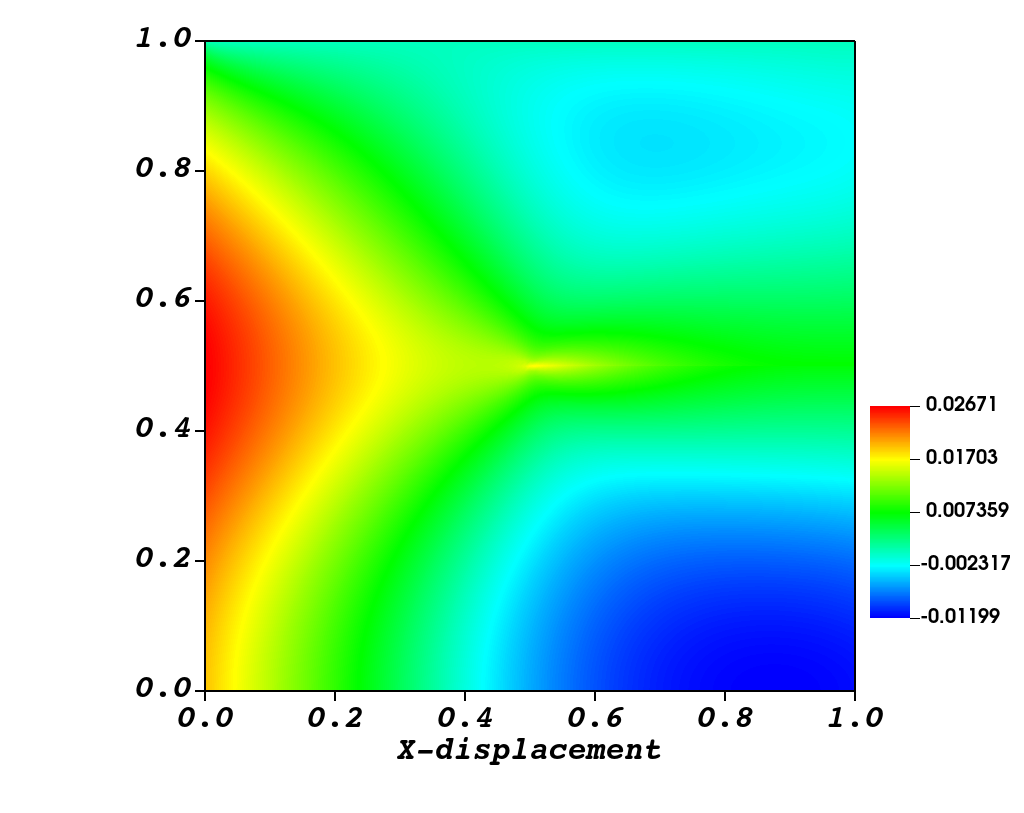}\quad
\includegraphics[scale=0.2]{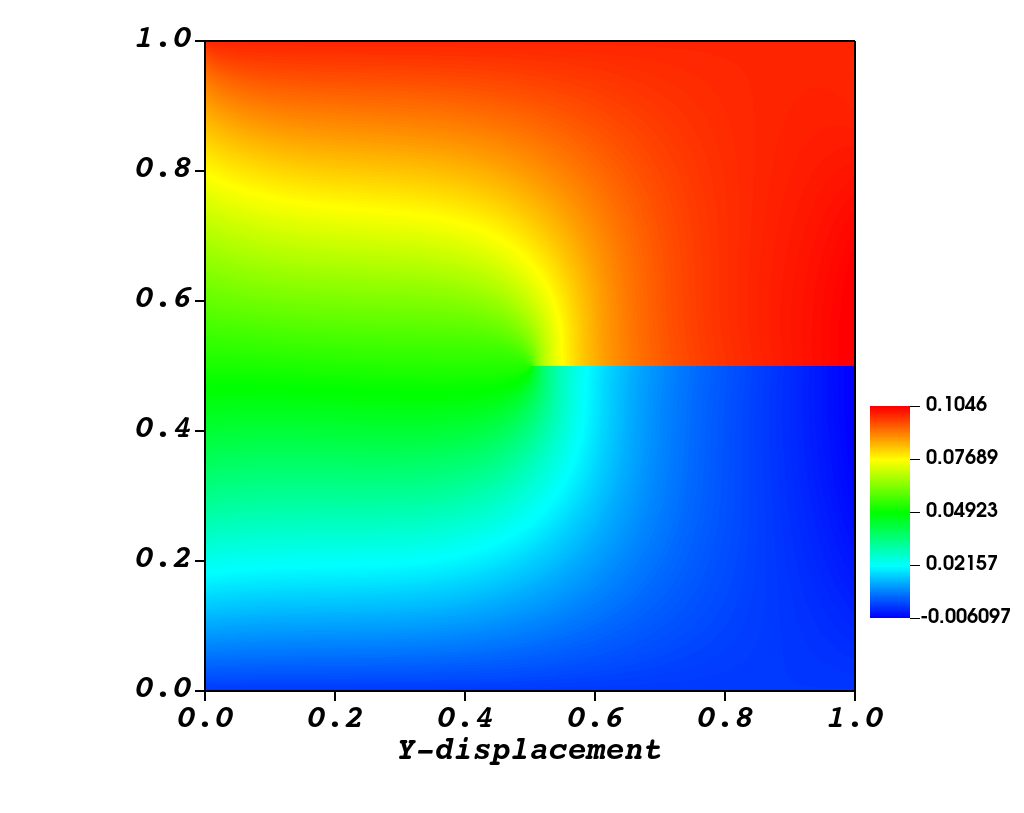}
\caption{horizontal and vertical displacements for $\theta=100,a=0.5,b=0.02$}\label{fig:dis_m1_const}
\end{figure}

Figure (\ref{fig:dis_m1_const}) illustrates the horizontal ($u_x$) and vertical ($u_y$) displacement fields resulting from the application of a constant temperature boundary condition. The mechanical boundary conditions are defined as follows:

\begin{itemize}
    \item \textbf{Top Boundary}: This boundary is constrained against horizontal movement, corresponding to a {homogeneous} Dirichlet boundary condition for horizontal displacement ($u_x=0$). However, it is allowed to move vertically, subject to a {non-homogeneous} condition.
    \item \textbf{Bottom Boundary}: Conversely, the bottom boundary is fixed to prevent vertical movement, representing a {homogeneous} Dirichlet boundary condition for vertical displacement ($u_y=0$).
\end{itemize}

This specific set of constraints is designed to simulate a scenario where the material is fixed at its base and restricted from expanding or contracting horizontally at the top. This allows for a clear analysis of thermally induced deformation.

\begin{figure}[H]
\centering 
\centering{}\includegraphics[scale=0.2]{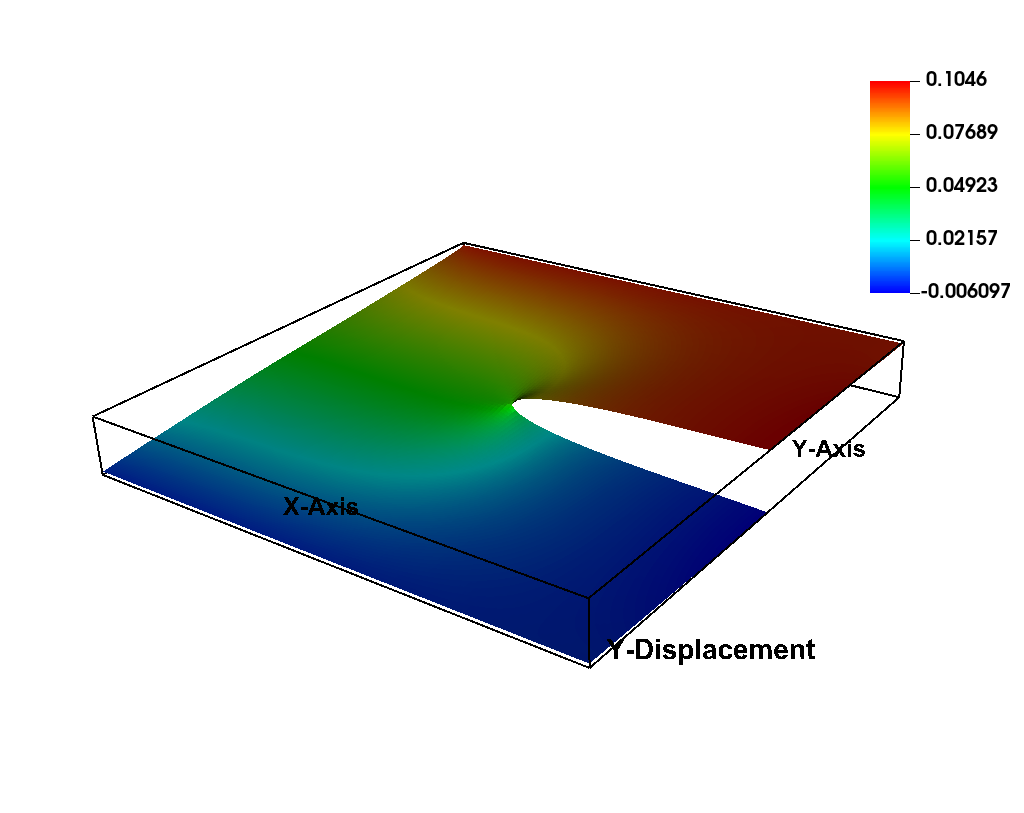}
\caption{Elevated $u_{y}$, $\theta=100,a=0.5,b=0.02$}\label{fig:y3d_m1_const}
\end{figure}

Figure~\ref{fig:y3d_m1_const} provides a three-dimensional visualization of the vertical displacement field, with particular emphasis on the deformation along the crack faces. The plot clearly shows that the most elevated vertical displacements are localized on the surfaces of the crack itself, which physically represents the crack opening under the applied thermo-mechanical load. A key observation from this profile is its distinct elliptical pattern. This shape is not arbitrary; it is the classic crack opening profile predicted by LEFM for a Mode I (opening mode) fracture. The emergence of this elliptical shape provides strong validation for the accuracy of the numerical model in capturing the fundamental mechanics of the fracture process.

\subsubsection{Parabolic Temperature}

In this section, we analyze the system's response to a non-uniform thermal load. Specifically, a {parabolic temperature profile} is prescribed along the bottom boundary using the Dirichlet condition $\theta(x) = 400x(1-x)$.

\begin{figure}[H]
\centering
\includegraphics[scale=0.2]{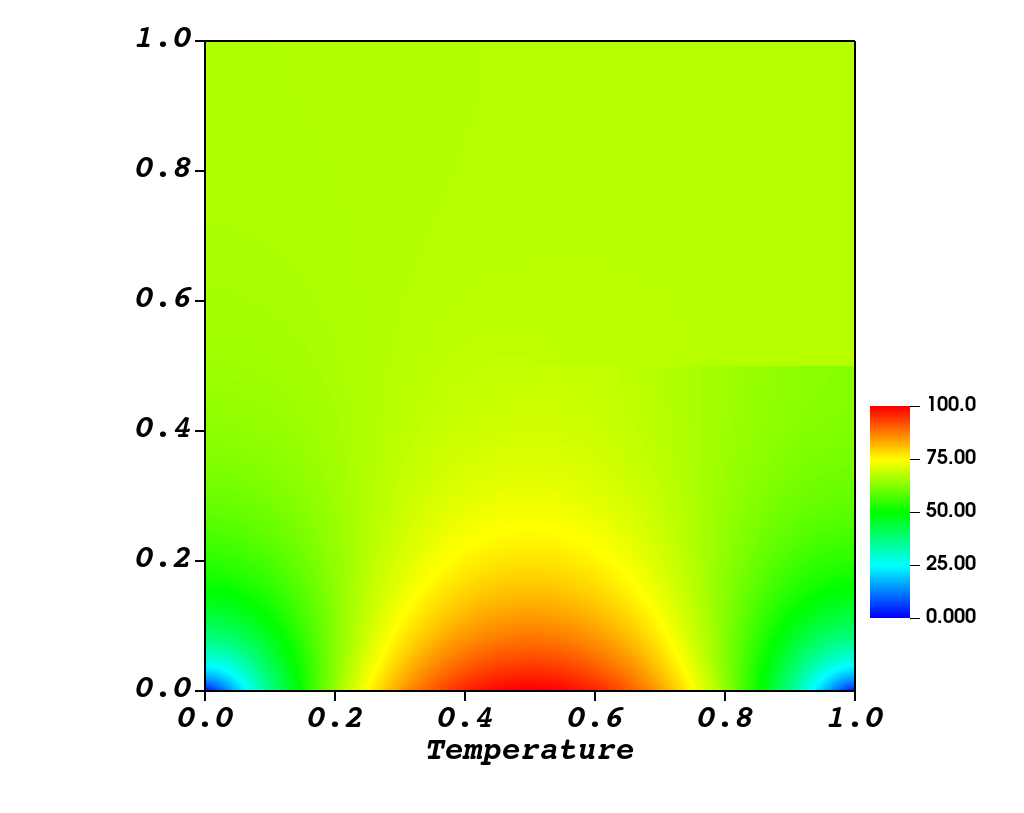}
\caption{Temperature profile for $\theta=400x(1-x),a=0.5,b=0.02$}\label{fig:Temp_m1_para}
\end{figure}

As illustrated in Figure~(\ref{fig:Temp_m1_para}), this function creates a distinct thermal gradient. The temperature is zero at the ends of the boundary ($x=0$ and $x=1$) and reaches its {maximum value of $\theta=100$} at the midpoint ($x=0.5$). This setup simulates a localized heating scenario, allowing for the investigation of thermally induced stresses originating from a non-homogeneous temperature field.

\begin{figure}[H]
\centering 
\includegraphics[scale=0.3]{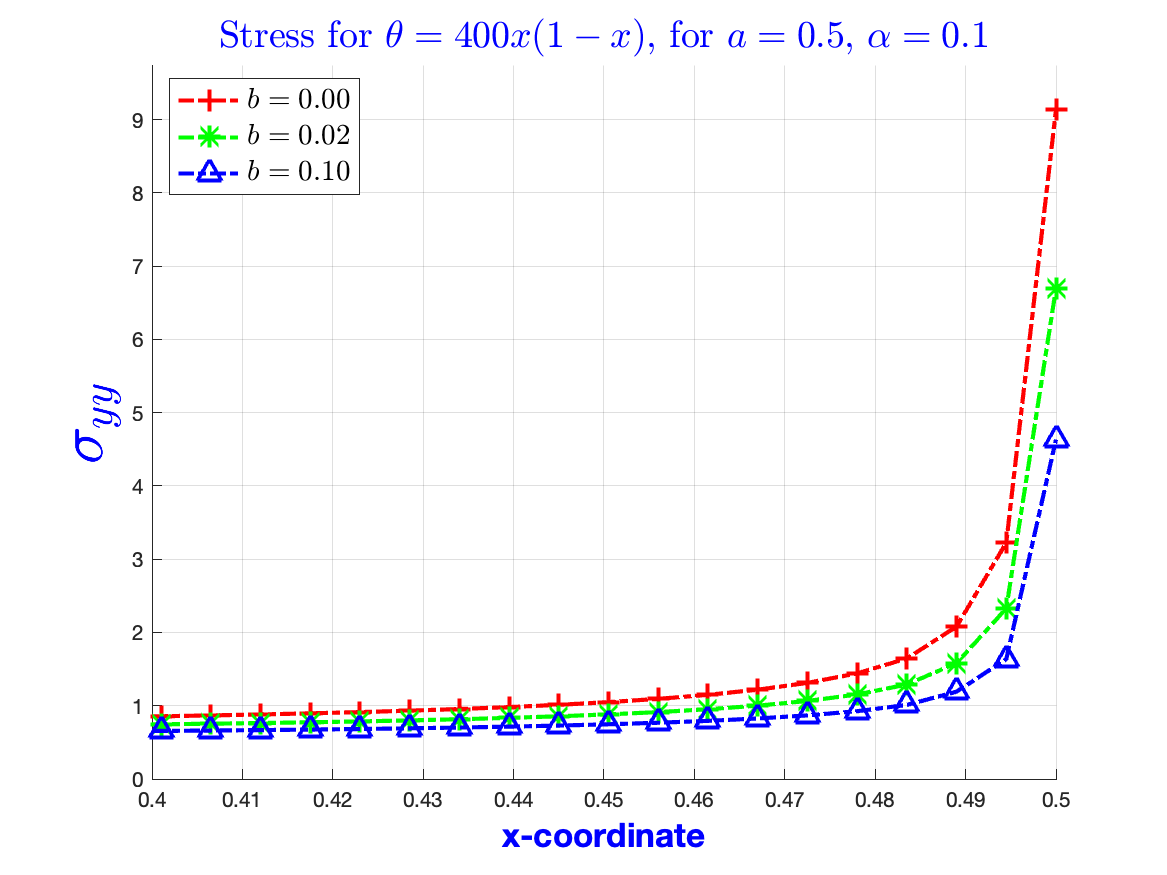}\quad
\includegraphics[scale=0.3]{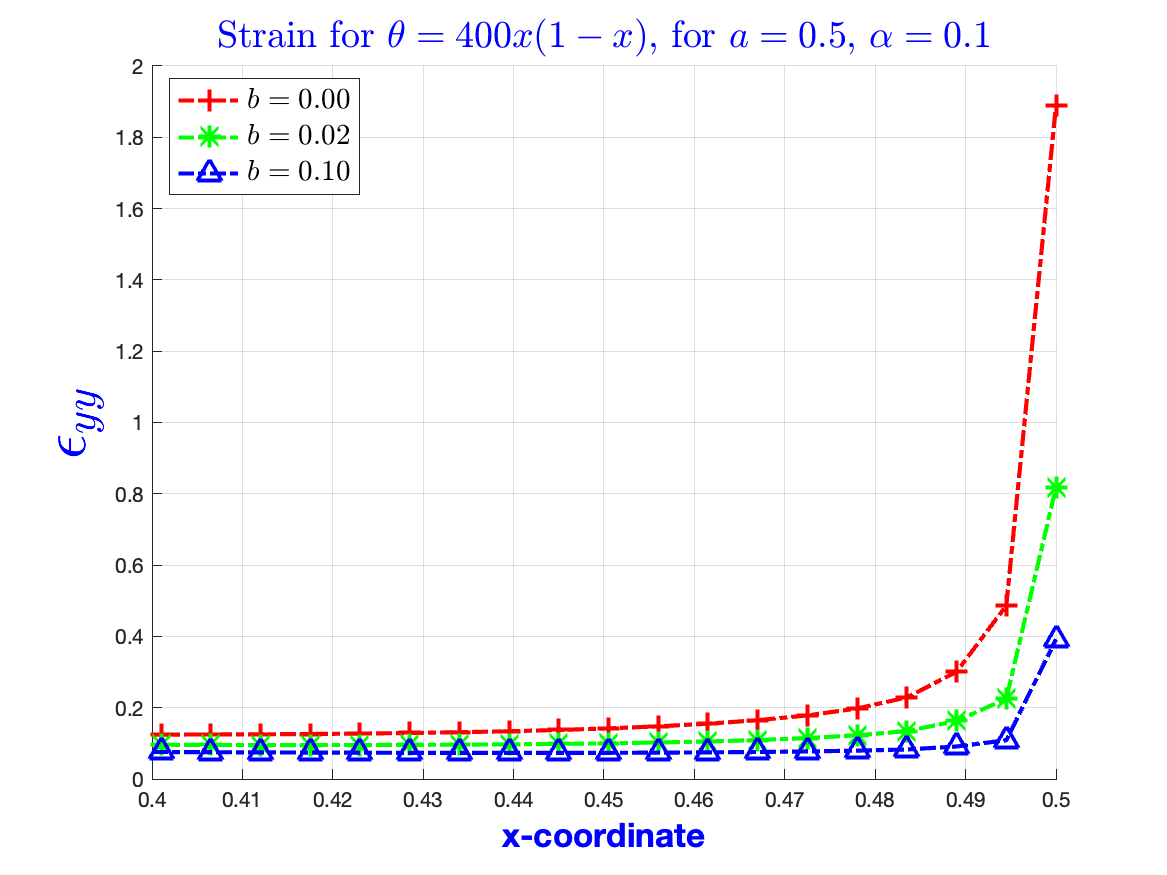}
\caption{Stress and strain for $b$ variation, $\theta=400x(1-x),a=0.5$}\label{fig:SS_m1_bv_para}
\end{figure}

Figure~\ref{fig:SS_m1_bv_para} illustrates the sensitivity of the stress and strain fields to the nonlinear material parameter $b$ under a parabolic temperature distribution applied to the bottom boundary. A clear inverse relationship is observed: as the value of $b$ increases, the magnitudes of both stress and strain decrease significantly. In the baseline linear case, corresponding to $b=0$, the stress and strain attain their maximum observed values of $9.93483$ and $1.96505$, respectively. This trend strongly suggests that higher values of the parameter $b$ enhance the material's resistance to fracture initiation by reducing local stress and strain concentrations.

\begin{figure}[H]
\centering 
\includegraphics[scale=0.3]{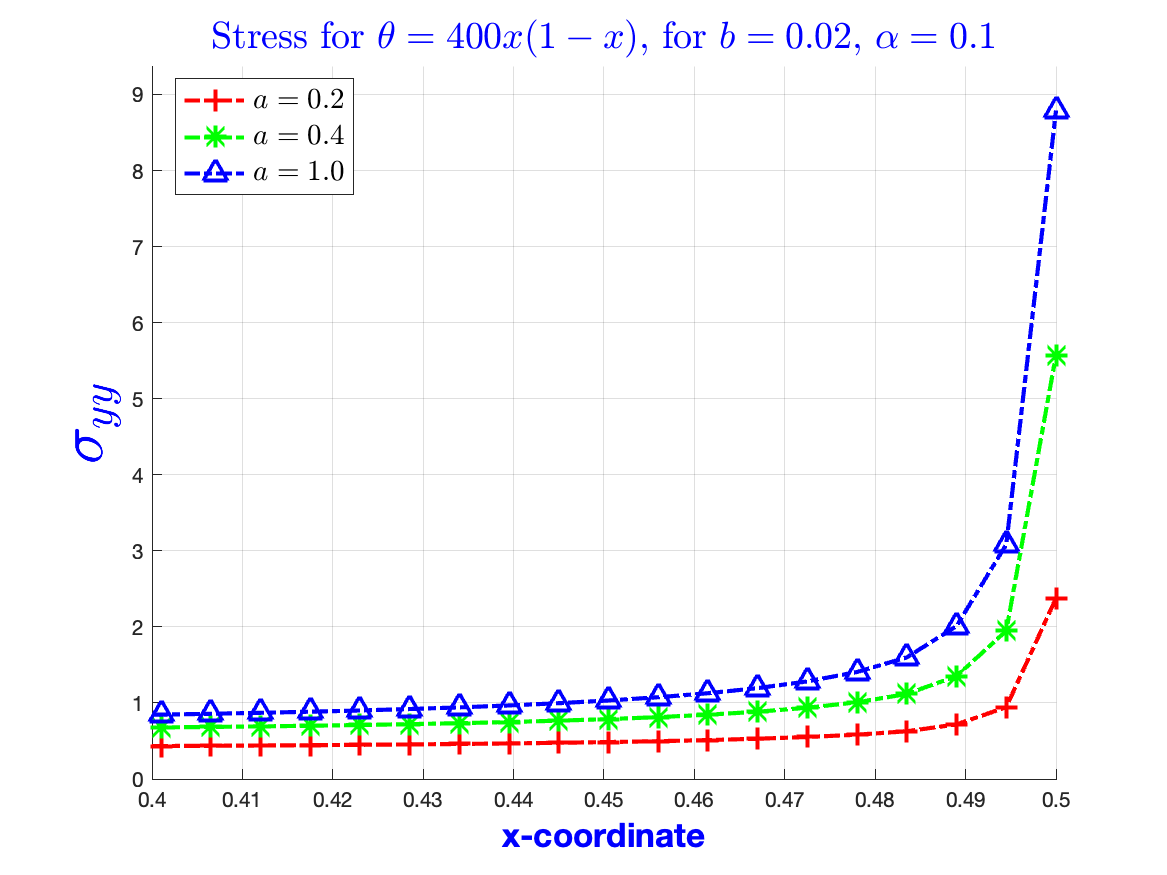}\quad
\includegraphics[scale=0.3]{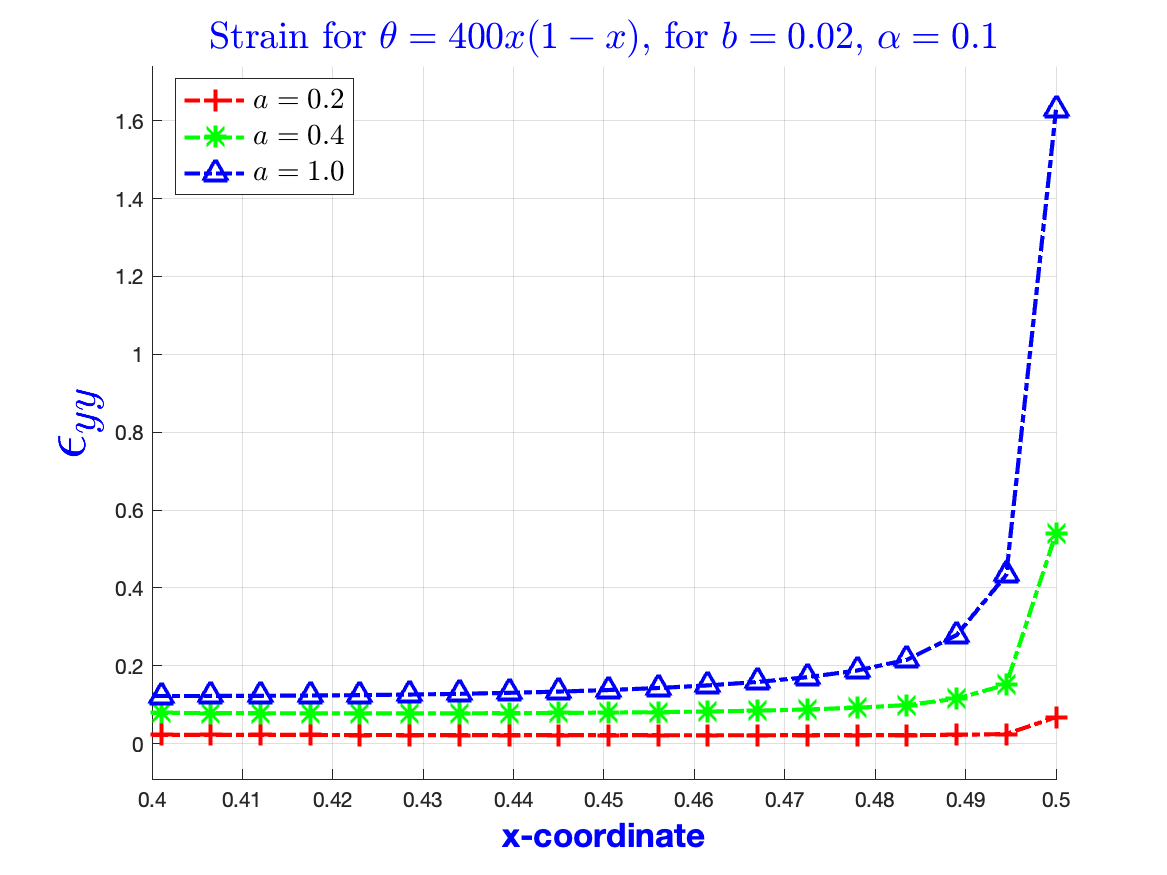}
\caption{Stress and strain for $a$ variation , $\theta=400x(1-x),b=0.02$}\label{fig:SS_m1_av_para}
\end{figure}

A contrasting behavior is presented in Figure~\ref{fig:SS_m1_av_para}, which illustrates the effect of varying the parameter $a$ while all other parameters are held fixed. A direct correlation is observed: higher values of $a$ result in a significant amplification of both the stress and strain. This indicates that the material becomes more vulnerable to damage, as the resulting higher stress concentrations mean that failure can be initiated by a comparatively lower applied force.

\begin{figure}[H]
\centering 
\includegraphics[scale=0.2]{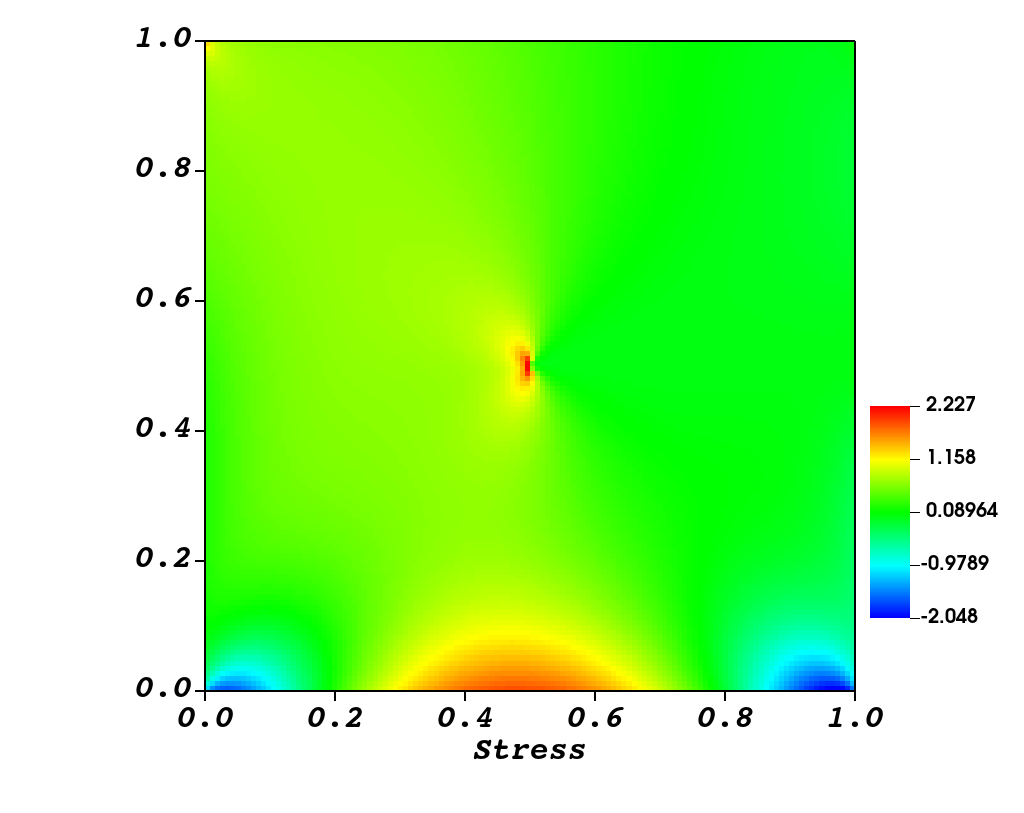}\quad
\includegraphics[scale=0.2]{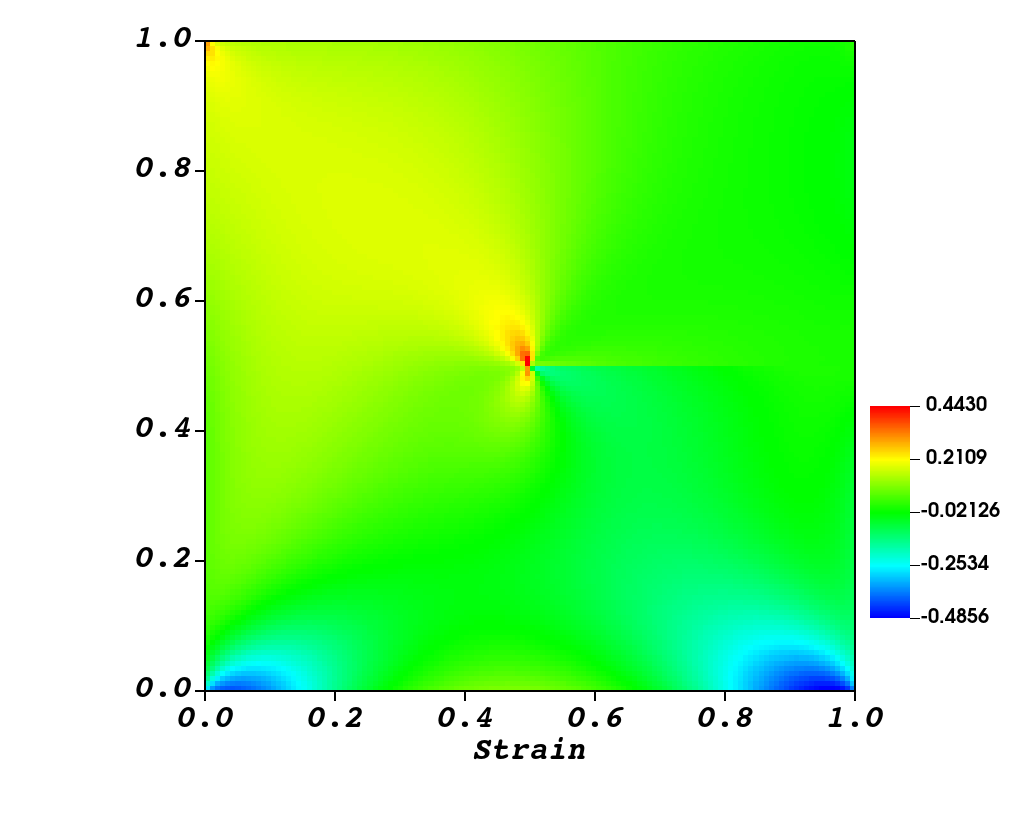}
\caption{Stress and Strain for $\theta=400x(1-x),a=0.5,b=0.02$}\label{fig:SS_m1_visit_para}
\end{figure}

Under thermo-mechanical loading, as illustrated in Figure~\ref{fig:SS_m1_visit_para}, both the stress and strain fields exhibit their maximum concentrations at the crack tip, as is expected in fracture analysis. A secondary, localized stress concentration is also observed near the middle of the bottom boundary. This effect is a direct consequence of the applied parabolic temperature profile, $\theta(x) = 400x(1-x)$, which induces significant thermal stresses in that region. Quantitatively, the analysis shows that the peak tensile stress reaches a value of $2.227$, while the maximum compressive stress is $-2.048$. Correspondingly, the maximum and minimum principal strain values are recorded as $0.4430$ and $-0.4856$, respectively.

\begin{figure}[H]
\centering 
\includegraphics[scale=0.2]{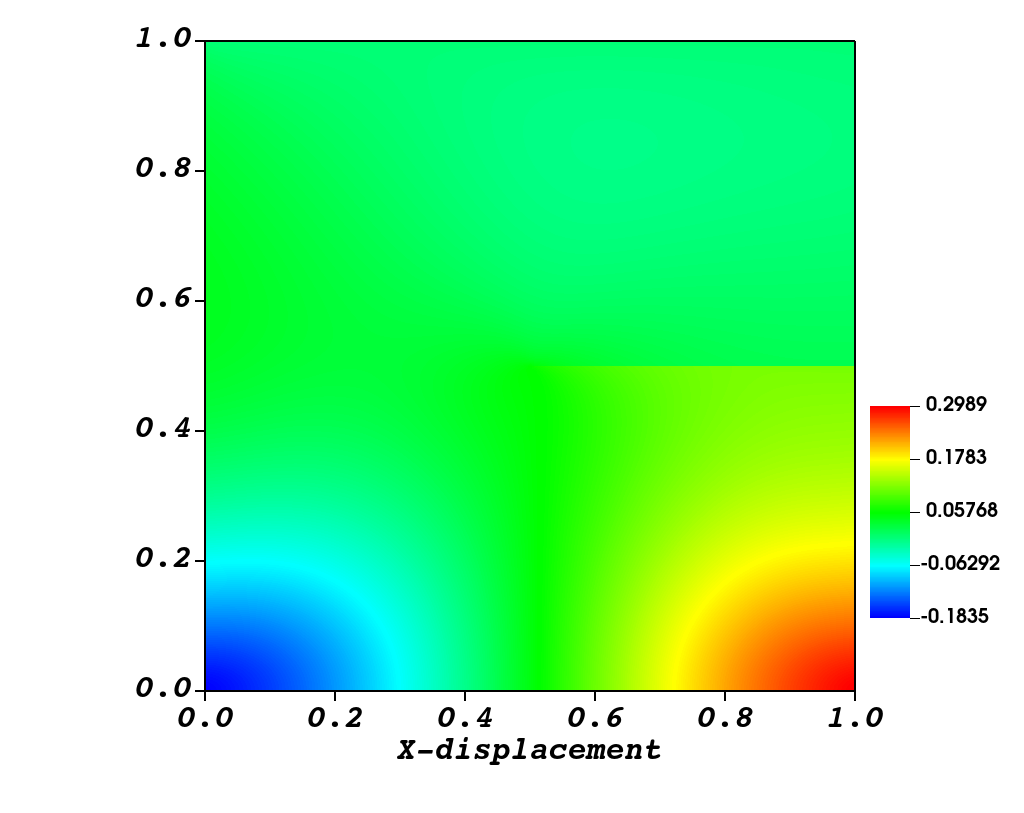}\quad
\includegraphics[scale=0.2]{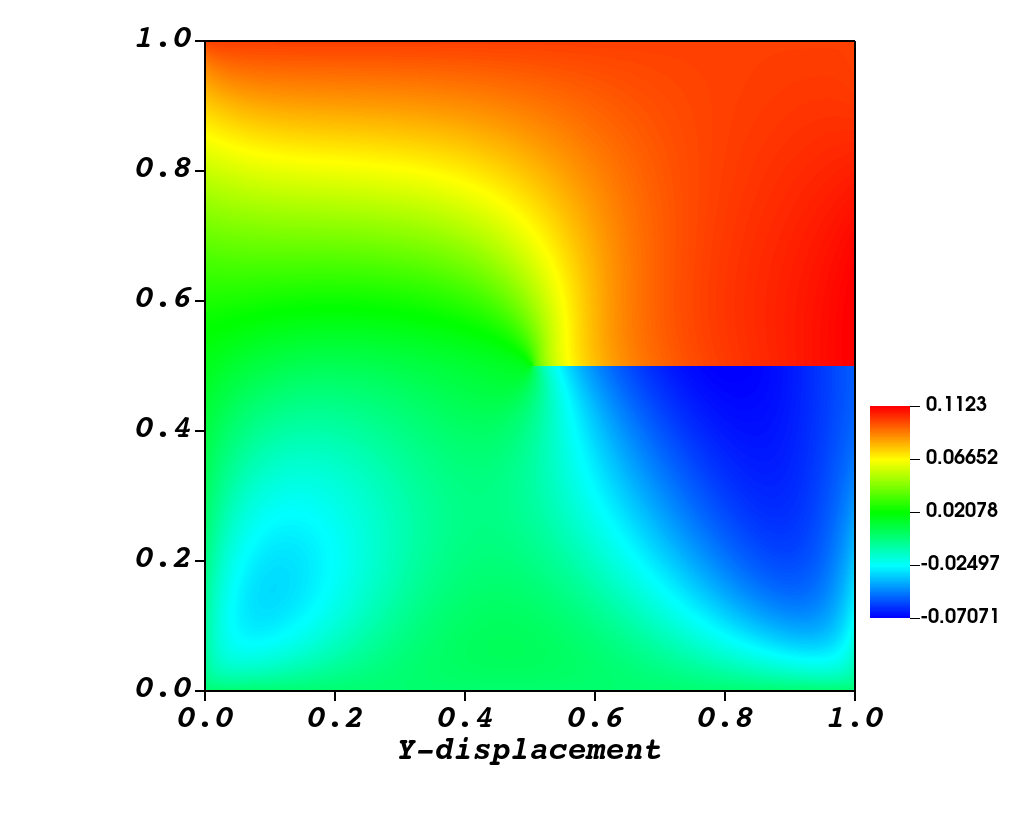}
\caption{horizontal and vertical displacements for $\theta=400x(1-x),a=0.5,b=0.02$}\label{fig:dis_m1_para}
\end{figure}

An analysis of the horizontal displacement component, shown in Figure~\ref{fig:dis_m1_para}, reveals a more complex pattern. Unlike the symmetric response typically seen in pure Mode-I scenarios, the horizontal field is clearly asymmetric. This asymmetry is a key feature of the solution, suggesting the activation of mixed-mode behavior at the crack tip.

\begin{figure}[H]
\centering
\includegraphics[scale=0.2]{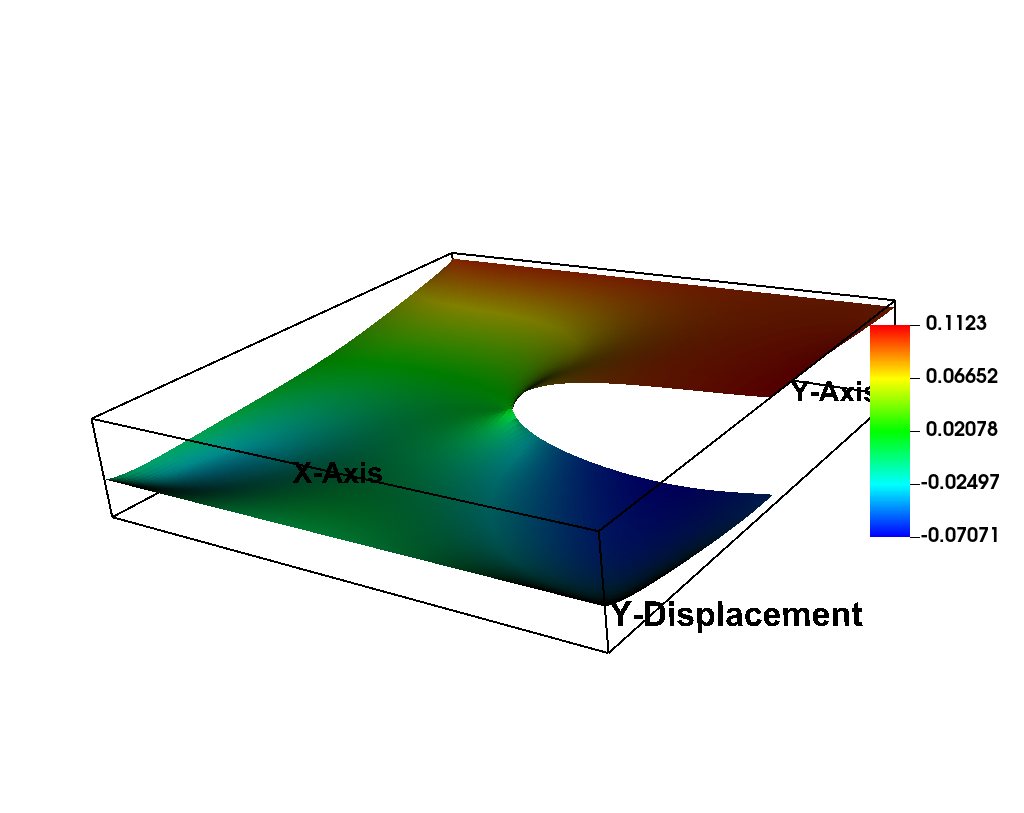}
\caption{$u_{y}$ elevated for $\theta=400x(1-x),a=0.5,b=0.02$ }\label{fig:y3d_m1_para}
\end{figure}

The comprehensive response of the system is illustrated in Figure~\ref{fig:y3d_m1_para}, which simultaneously displays the temperature distribution alongside the resulting horizontal and vertical displacement fields. The vertical displacement is of particular interest, as its profile along the crack surfaces clearly traces an elliptical opening. This specific geometry is the hallmark of a classic Mode I (opening mode) fracture, confirming the expected deformation behavior under the applied thermo-mechanical conditions.

\subsection{Fiber orientation orthogonal to crack}

This section analyzes the scenario where the material's anisotropic properties are oriented perpendicular to the crack plane. Specifically, the fiber reinforcement is aligned along the y-axis, which is defined as the axis of material symmetry. This preferred material direction is mathematically represented by the second-order structural tensor $\boldsymbol{M}=\boldsymbol{e}_{2}\otimes\boldsymbol{e}_{2}$, where $\boldsymbol{e}_{2}$ denotes the unit vector along the y-axis.

\subsubsection{Constant Temperature}

\begin{figure}[H]
\centering
 \includegraphics[scale=0.2]{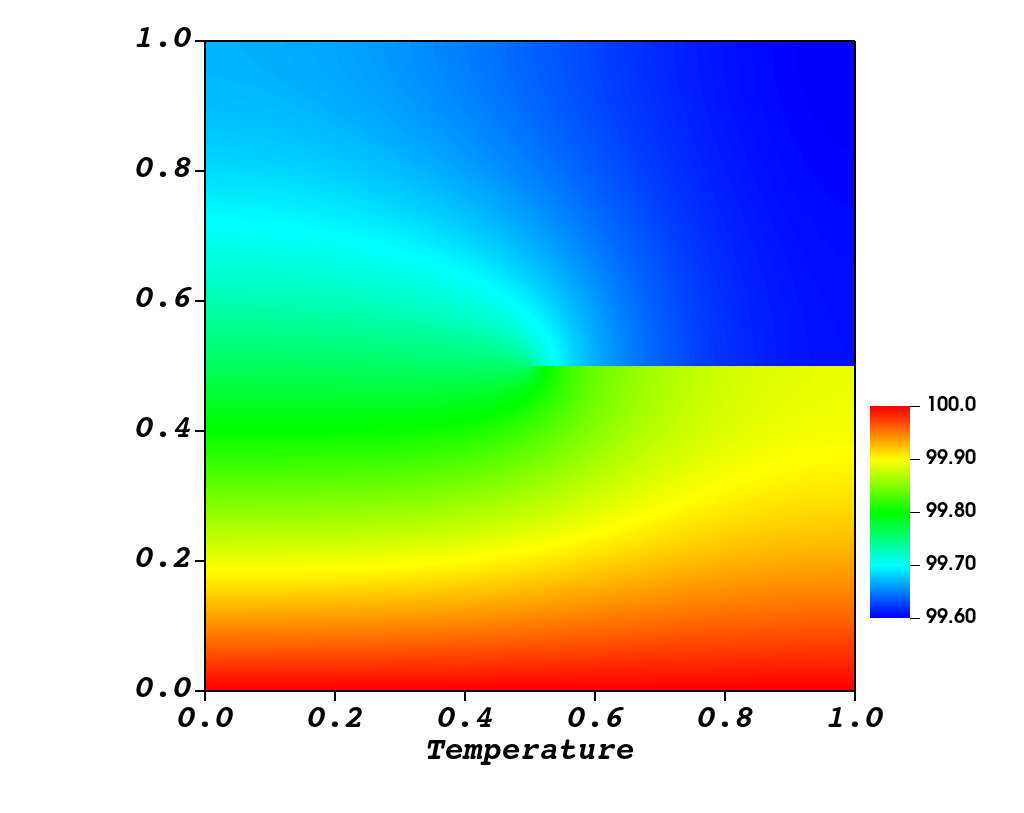}
\caption{Temperature profile, $\theta=100,a=0.5,b=0.02$}\label{fig:temp_m2_const}
\end{figure}

The temperature profile for the case of a prescribed constant temperature Dirichlet boundary condition on the bottom boundary is presented in Figure~\ref{fig:temp_m2_const}. As shown, the domain remains in a nearly isothermal state, exhibiting only a slight temperature variation. Quantitatively, the temperature ranges from a maximum of $100.00$ at the specified boundary to a minimum of $99.60$ elsewhere, confirming the minimal thermal gradient throughout the body.

\begin{figure}[H]
\centering 
\includegraphics[scale=0.3]{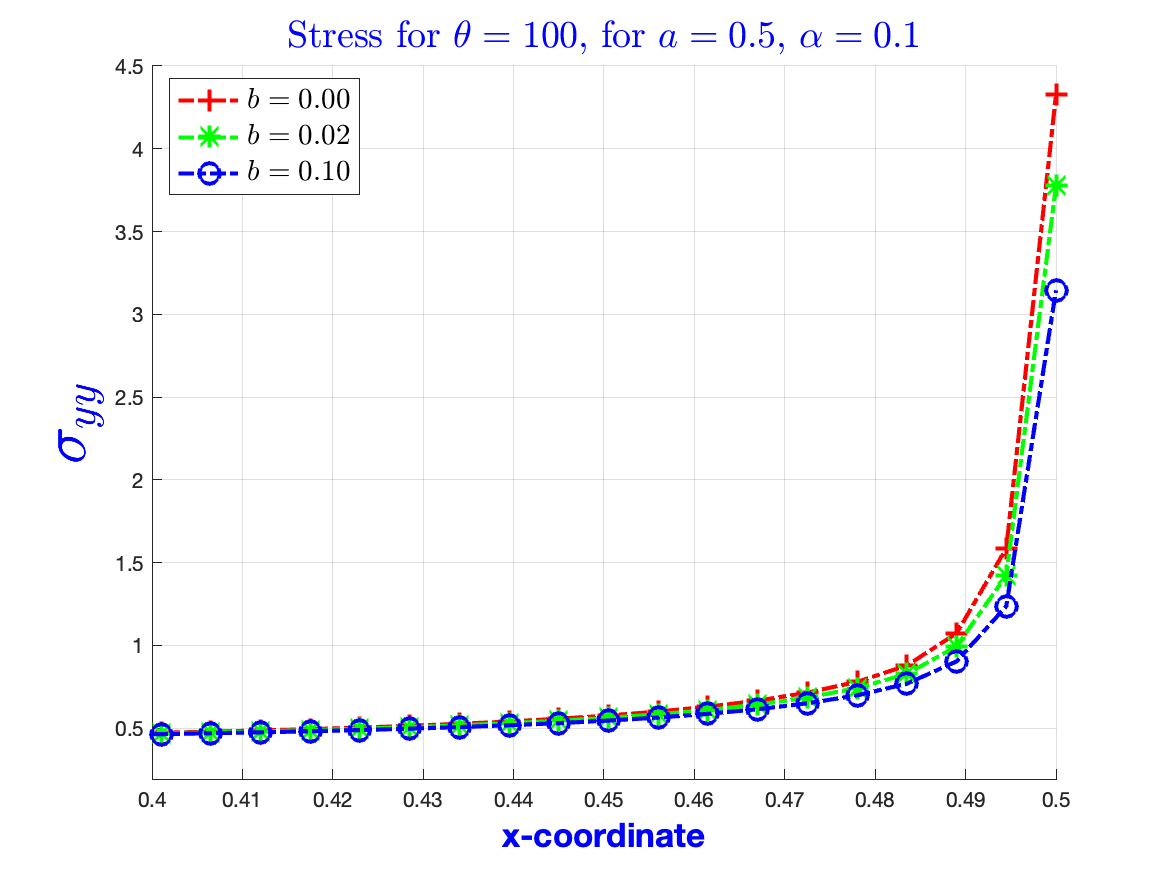}\quad
\includegraphics[scale=0.3]{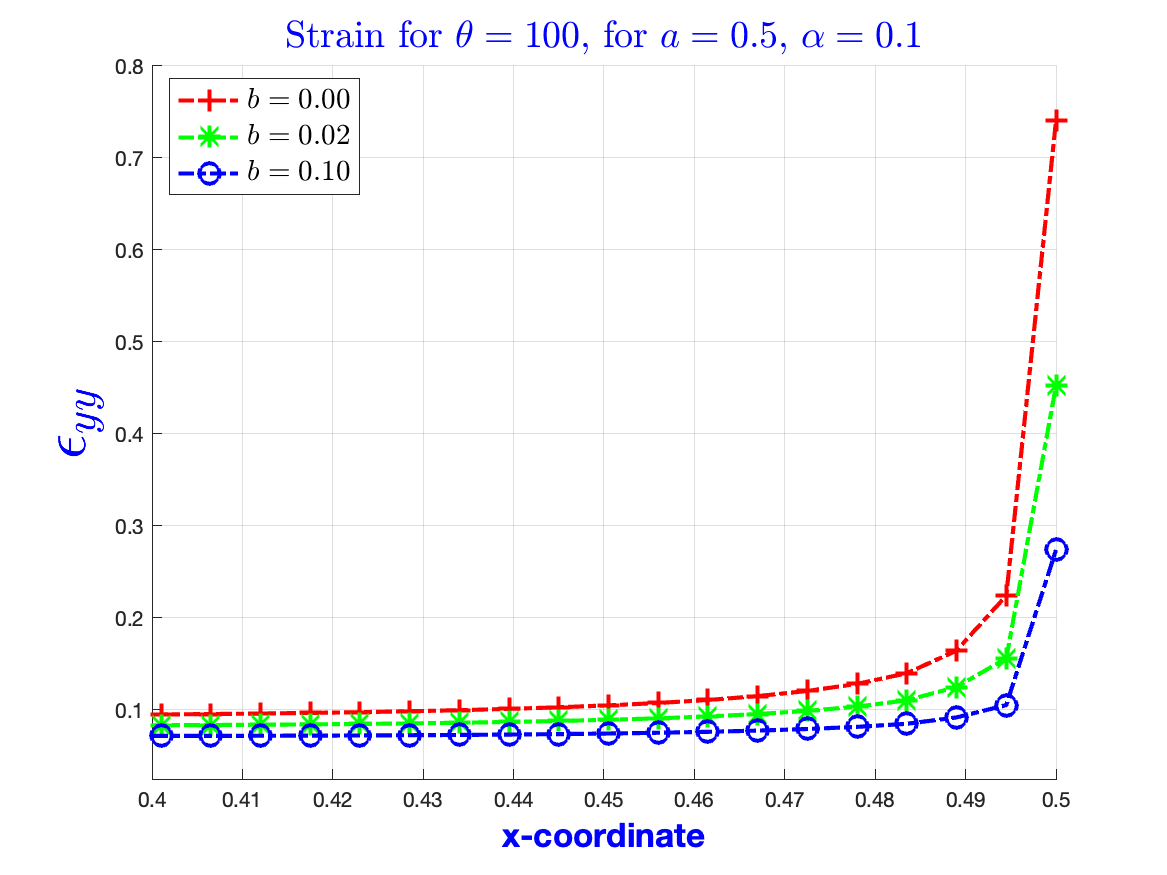}
\caption{Stress and strain for $b$ variation, $\theta=100,a=0.5$}\label{fig:SS_m2_bv_const}
\end{figure}

An analysis of the influence of the nonlinear parameter $b$ is presented in Figure~\ref{fig:SS_m2_bv_const}, revealing a clear trend of mitigation. Specifically, an inverse correlation exists: increasing the value of $b$ leads to a systematic reduction in both the stress and strain concentrations. The peak values for stress ($4.69695$) and strain ($0.769262$) are recorded for the reference linear elastic case, where $b=0$.

\begin{figure}[H]
\centering 
\includegraphics[scale=0.3]{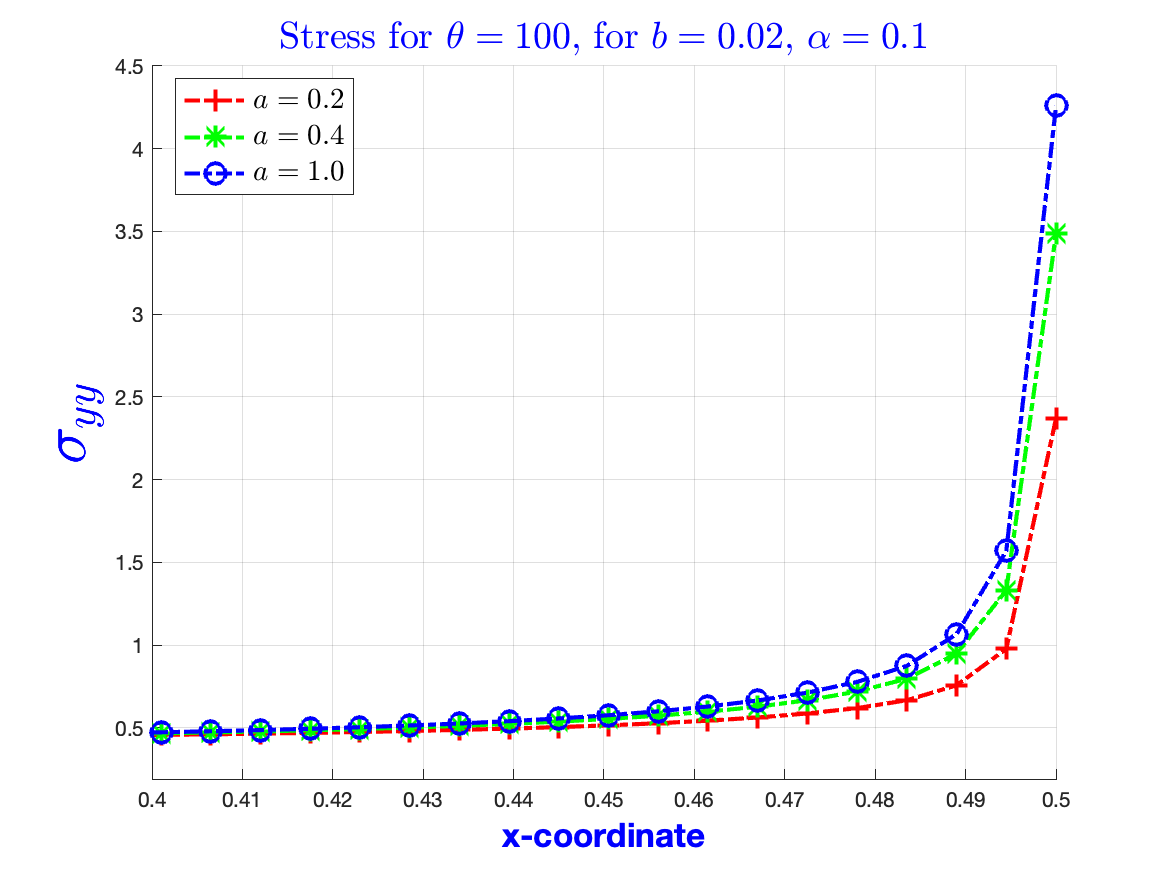}\quad
\includegraphics[scale=0.3]{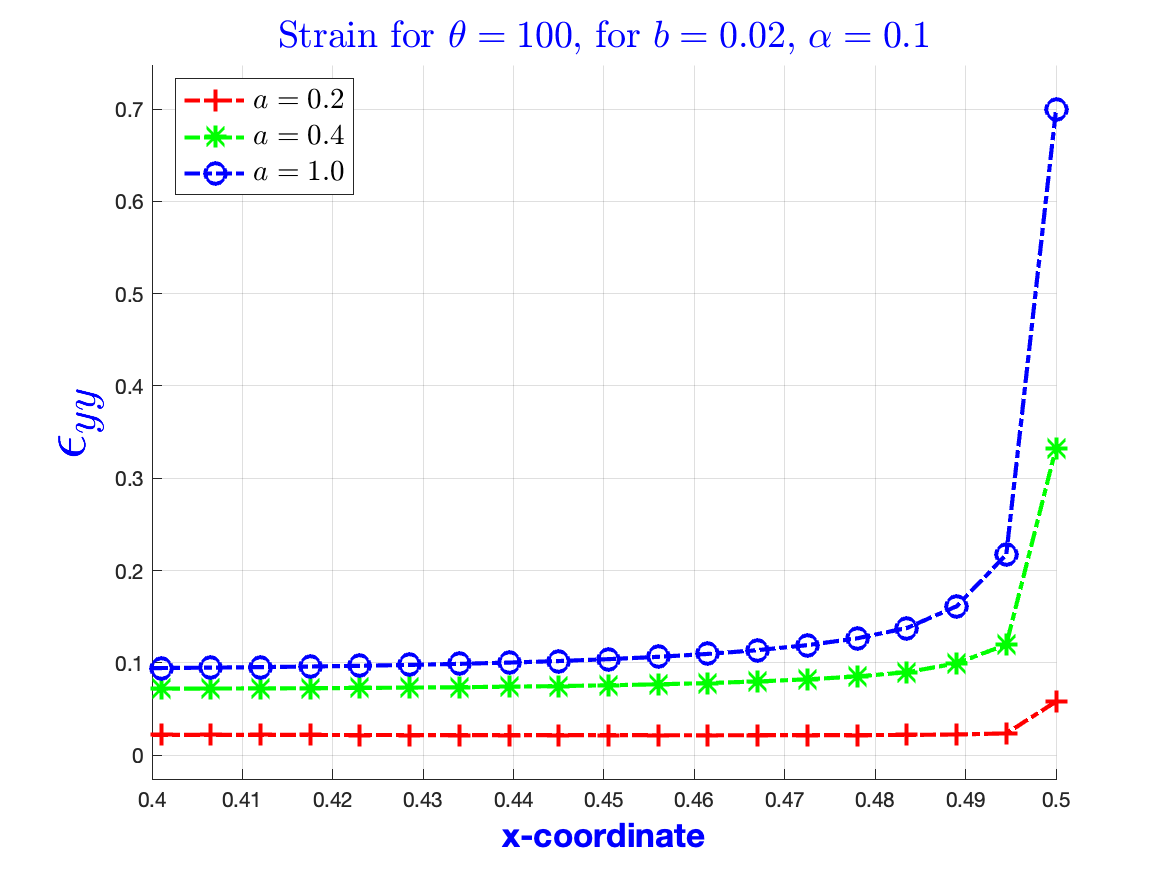}
\caption{Stress and strain for $a$ variation, $\theta=100,b=0.02$}\label{fig:SS_m2_av_const}
\end{figure}

In contrast, an analysis of the parameter $a$ reveals an opposing trend, as illustrated in Figure~\ref{fig:SS_m2_av_const}. A direct correlation is observed: higher values of $a$ lead to a greater mechanical response, amplifying both the stress and strain fields under the same constant temperature boundary condition. This behavior suggests that increasing $a$ exacerbates the material's susceptibility to fracture by elevating local stress concentrations.

\begin{figure}[H]
\centering 
\includegraphics[scale=0.2]{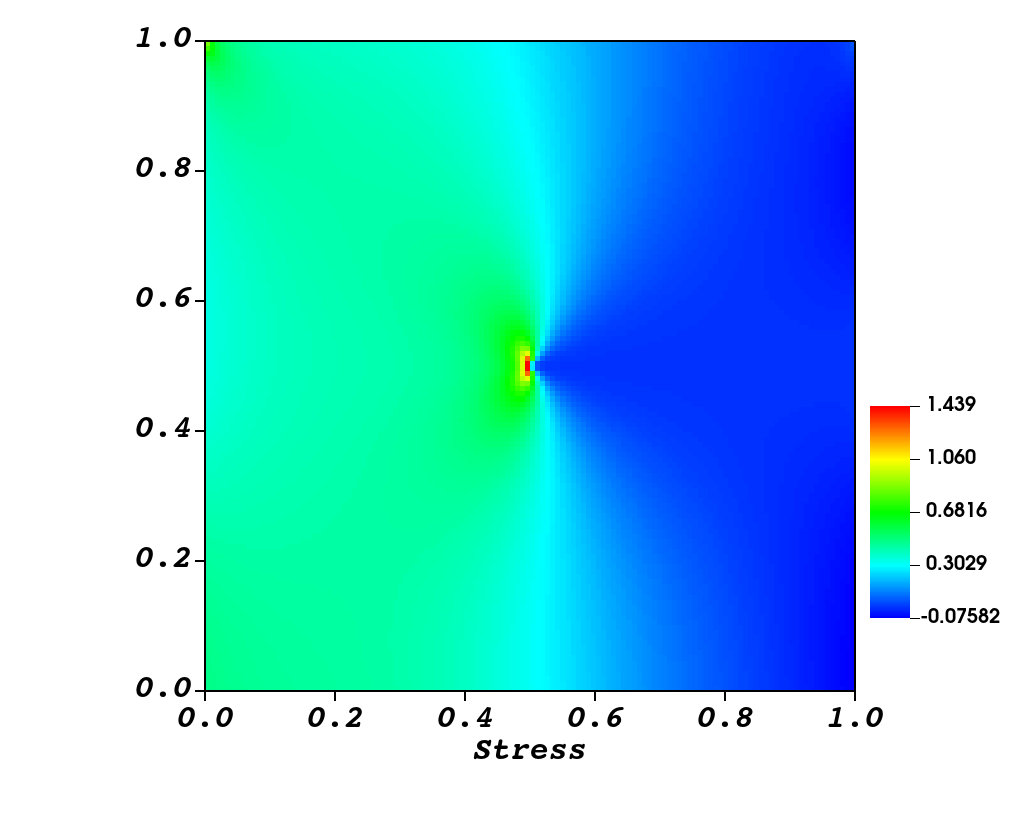}\quad
\includegraphics[scale=0.2]{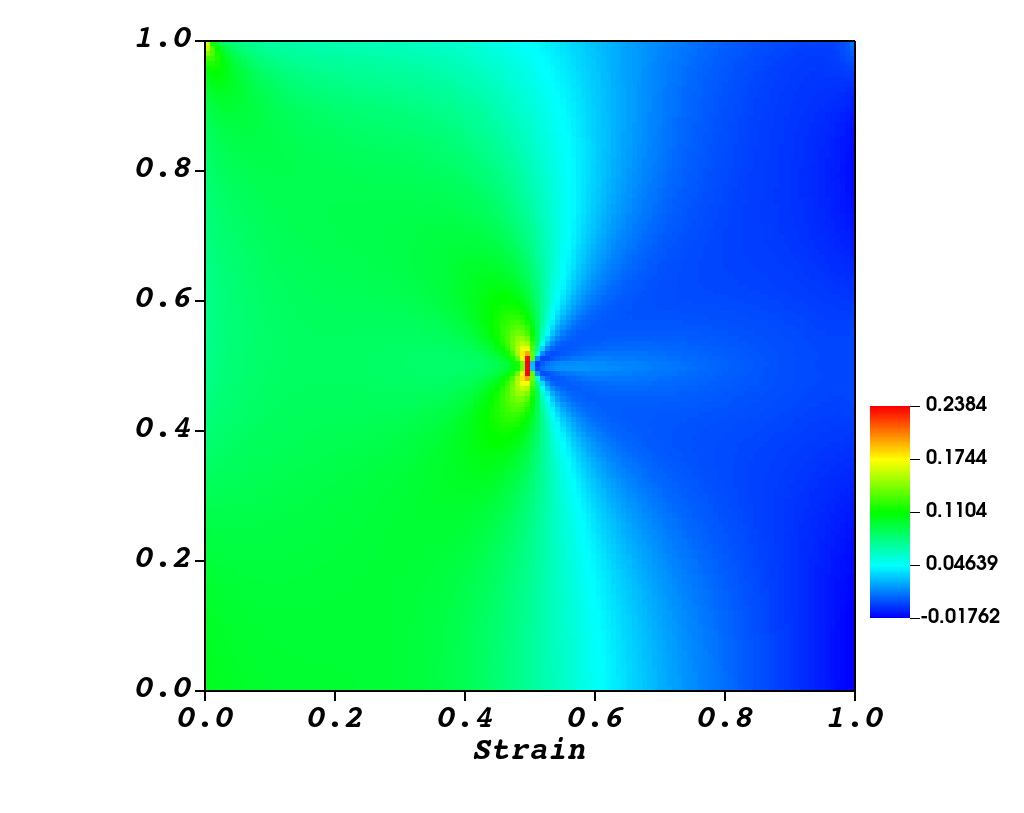}
\caption{Stress and strain for $\theta=100,a=0.5,b=0.02$}\label{fig:SS_m2_visit_const}
\end{figure}

The spatial distribution of the stress and strain fields under a constant temperature boundary condition is illustrated in Figure~\ref{fig:SS_m2_visit_const}. The results show that both mechanical fields are highly localized, reaching their peak values at the crack tip, which is characteristic of a stress singularity in fracture mechanics. Away from this region of high concentration, the fields decay rapidly. Notably, the stress and strain values are lowest in the far-field, particularly in the vicinity of the right-hand boundary of the domain.


\begin{figure}[H]
\centering 
\includegraphics[scale=0.2]{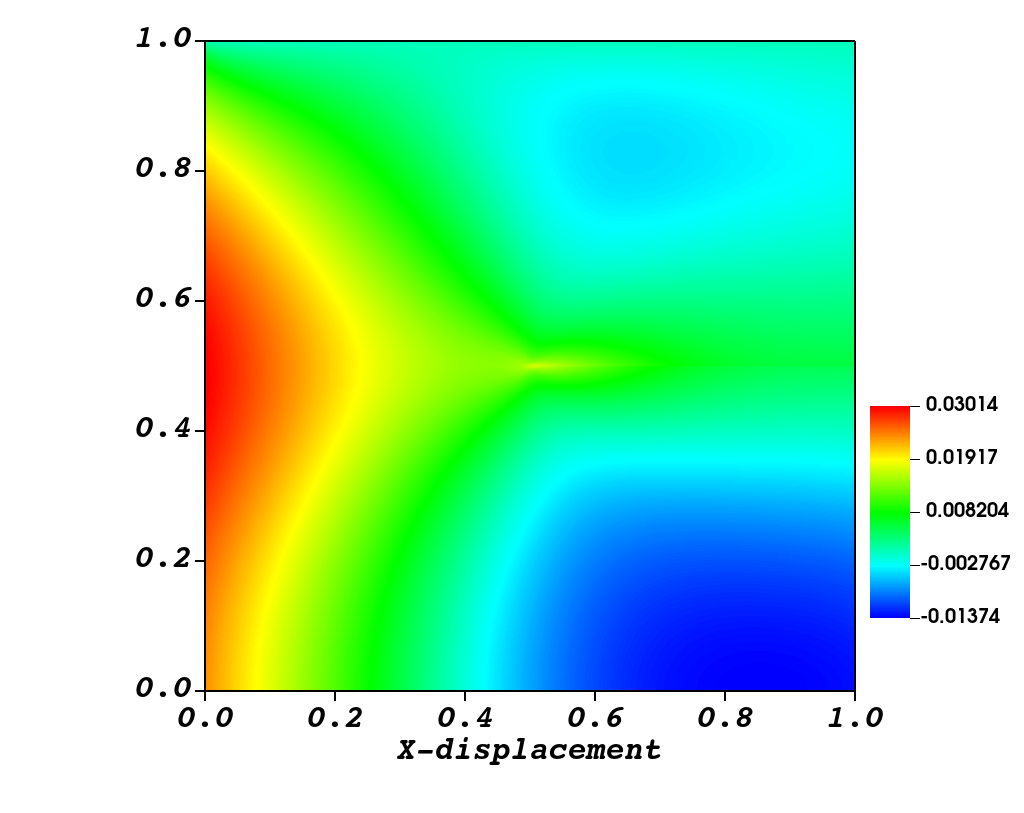}\quad
\includegraphics[scale=0.2]{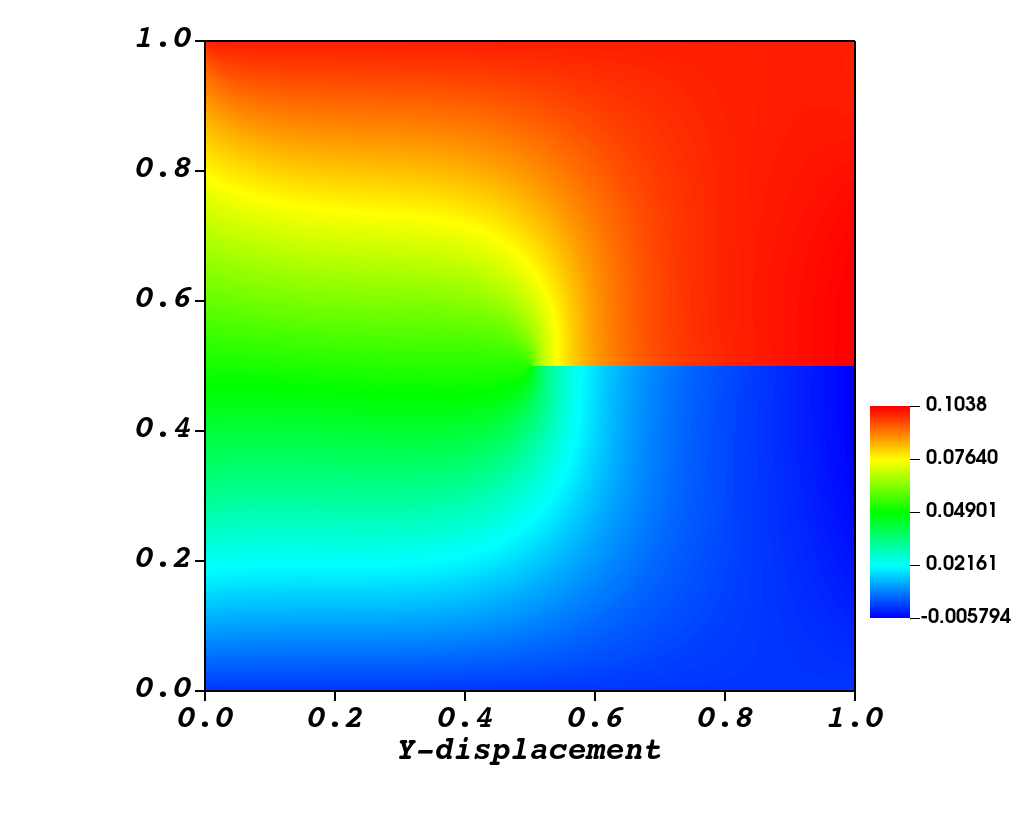}
\caption{Horizontal and vertical displacements for $\theta=100,a=0.5,b=0.02$}\label{fig:dis_m2_const}
\end{figure}

Figure~(\ref{fig:dis_m2_const}) suggests that the horizontal displacement 
is slightly higher near the left boundary and vertical displacement 
is higher near the upper half of the right boundary.

\begin{figure}[H]
\centering
\includegraphics[scale=0.2]{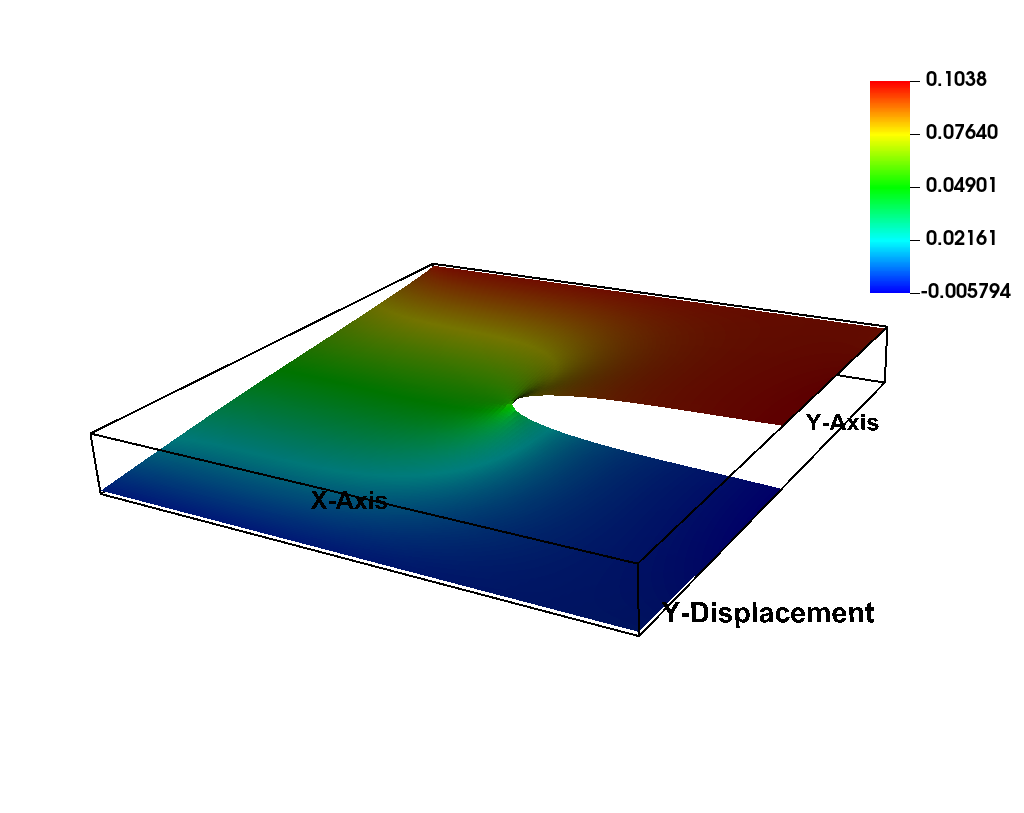} 
\caption{$u_{y}$ elevated, $\theta=100,a=0.5,b=0.02$}\label{fig:y3d_m2_const}
\end{figure}

As seen in Figure~(\ref{fig:y3d_m2_const}), the vertical displacement under a uniform bottom boundary is not symmetric. This asymmetry, evident in the crack profile and contour, is caused by the prescribed displacement boundary conditions.

\subsubsection{Parabolic Temperature}

\begin{figure}[H]
\centering
 \includegraphics[scale=0.2]{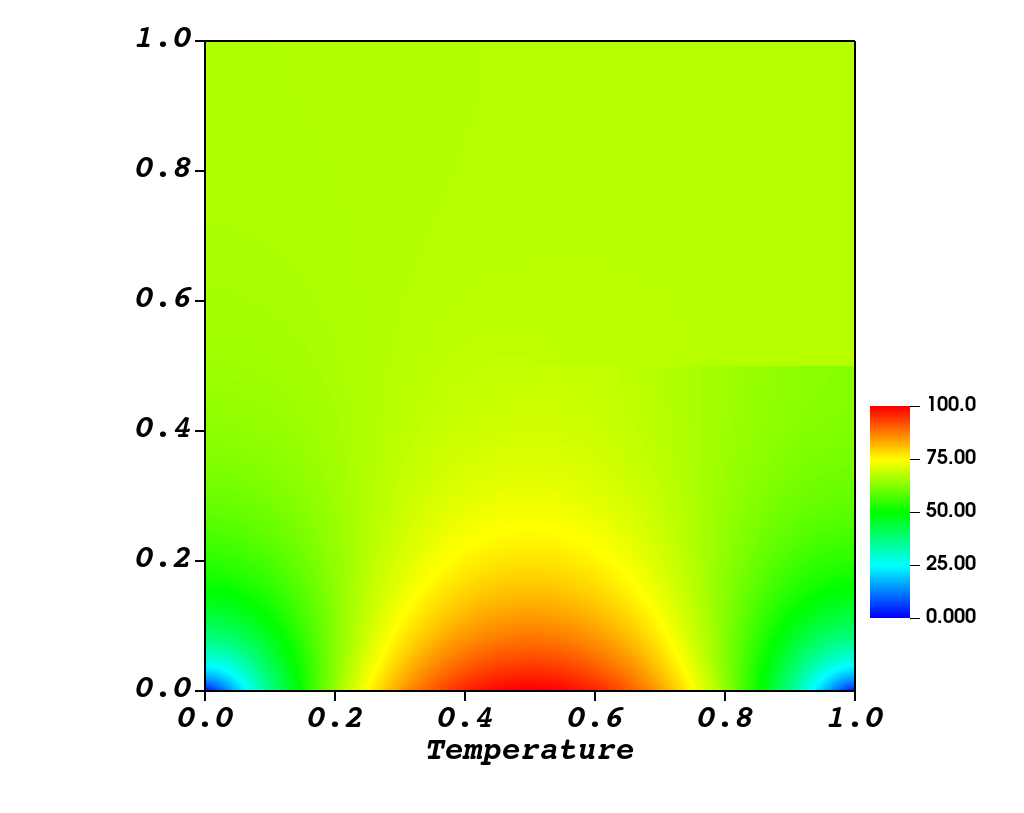}
\centering{}\caption{Temperature profile for $\theta=400x(1-x),a=0.5,b=0.02$}\label{fig:temp_m2_para}
\end{figure}

Figure~(\ref{fig:temp_m2_para}) illustrates the resulting temperature field, which peaks at the midpoint of the lower boundary ($x=0.5$) and decays progressively in the vertical direction. This observed thermal profile is characteristic of the applied parabolic Dirichlet boundary condition, which concentrates the maximum thermal energy at the center of the bottom edge.

\begin{figure}[H]
\centering
\includegraphics[scale=0.3]{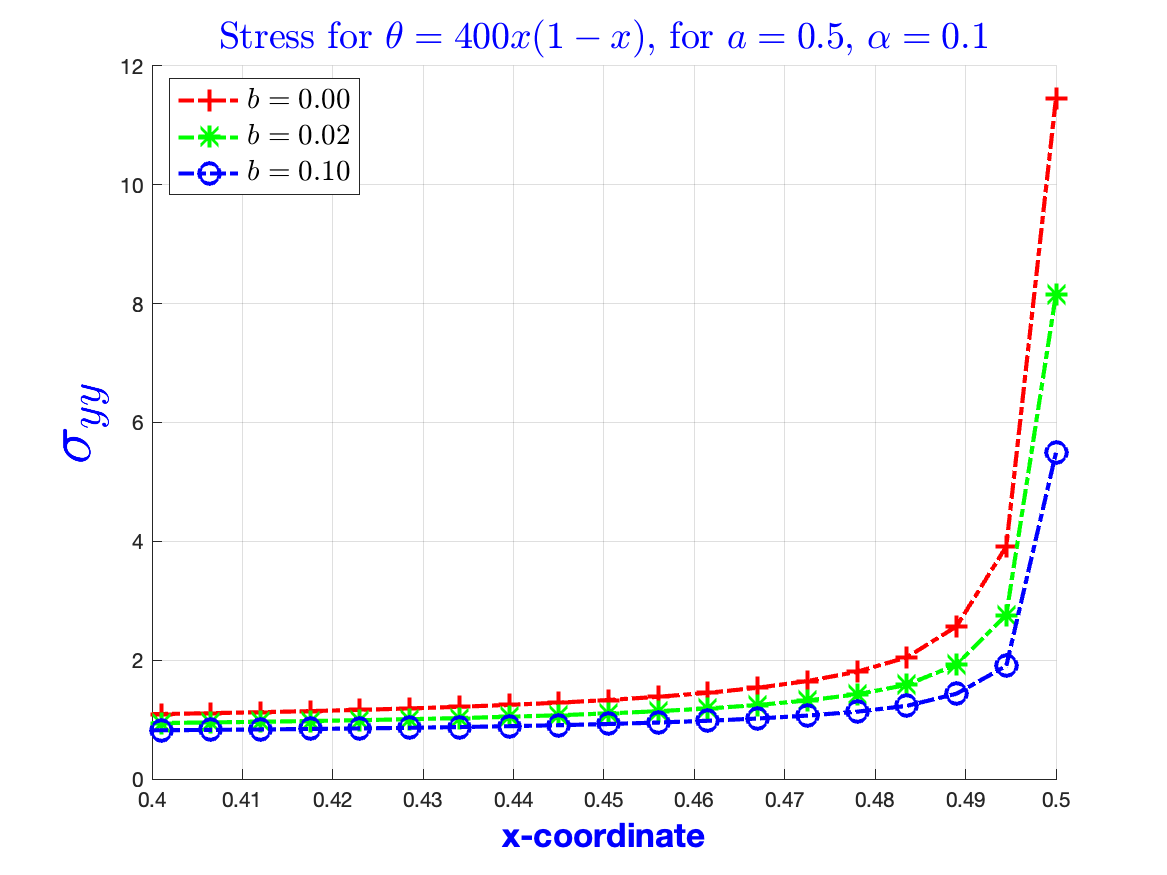}\quad
\includegraphics[scale=0.3]{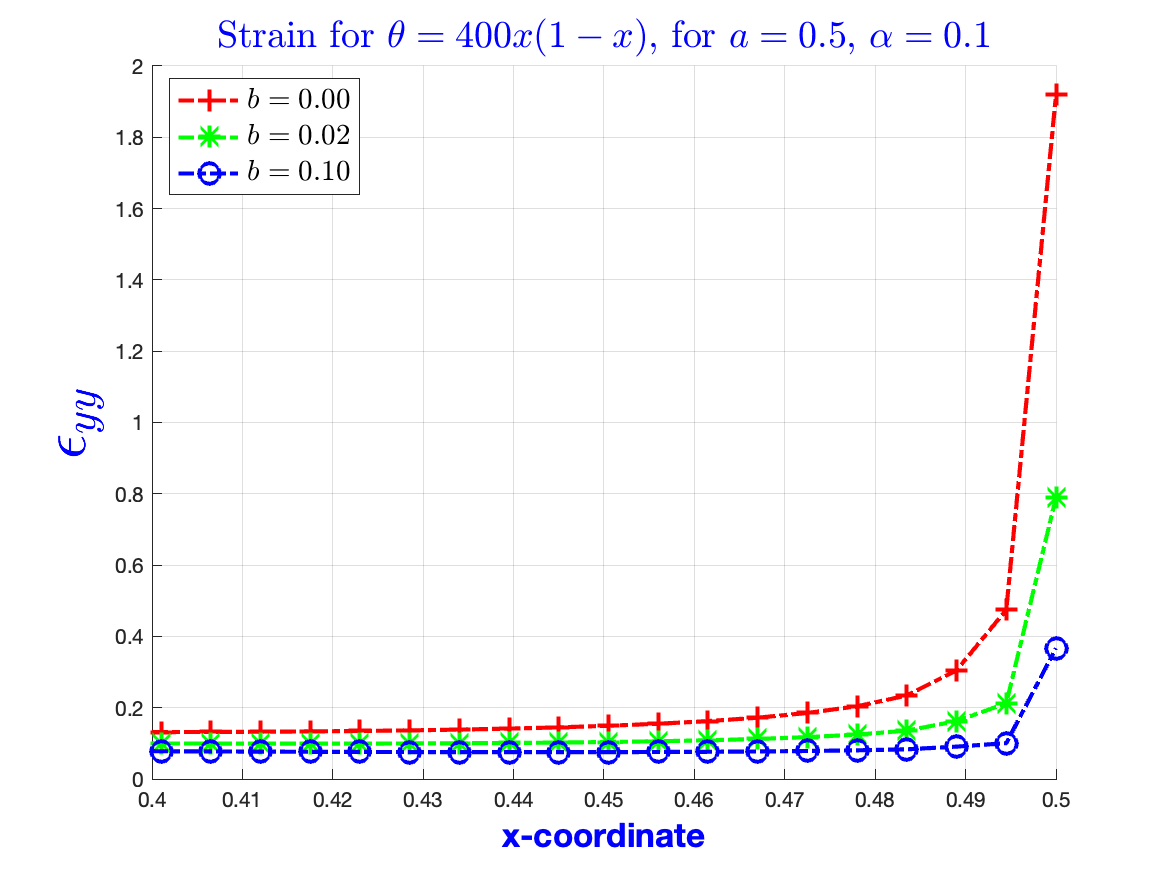}
\caption{Stress and strain for $b$ variation, $\theta=400x(1-x),a=0.5$}\label{fig:SS_m2_bv_para}
\end{figure}

Figure~(\ref{fig:SS_m2_bv_para}) demonstrates that the parameter $b$ is inversely correlated with the resulting stress and strain. This behavior suggests an improvement in the material's fracture toughness. Specifically, the lower stress levels imply that less strain energy is available to drive crack growth. Therefore, a material with a higher $b$ value is less susceptible to crack propagation under the applied loading conditions.

\begin{figure}[H]
\centering
\includegraphics[scale=0.3]{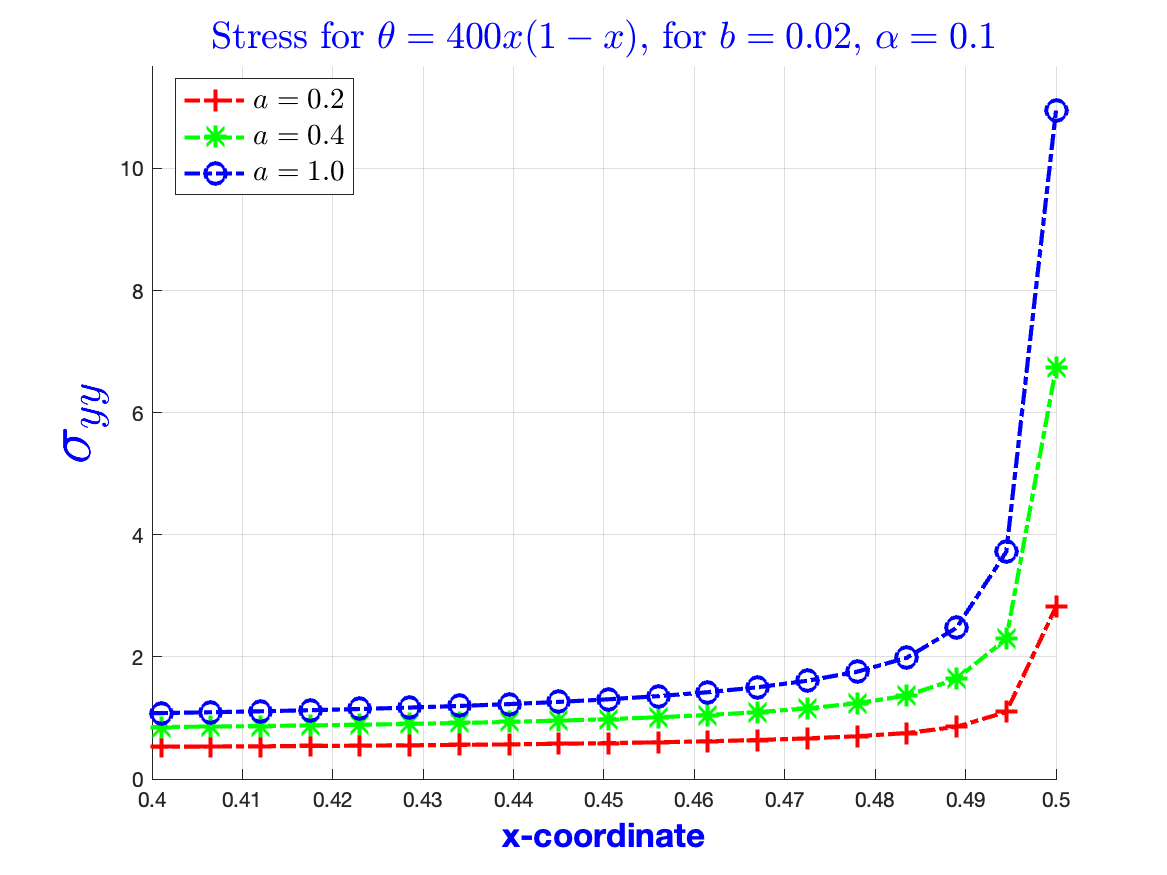}\quad
\includegraphics[scale=0.3]{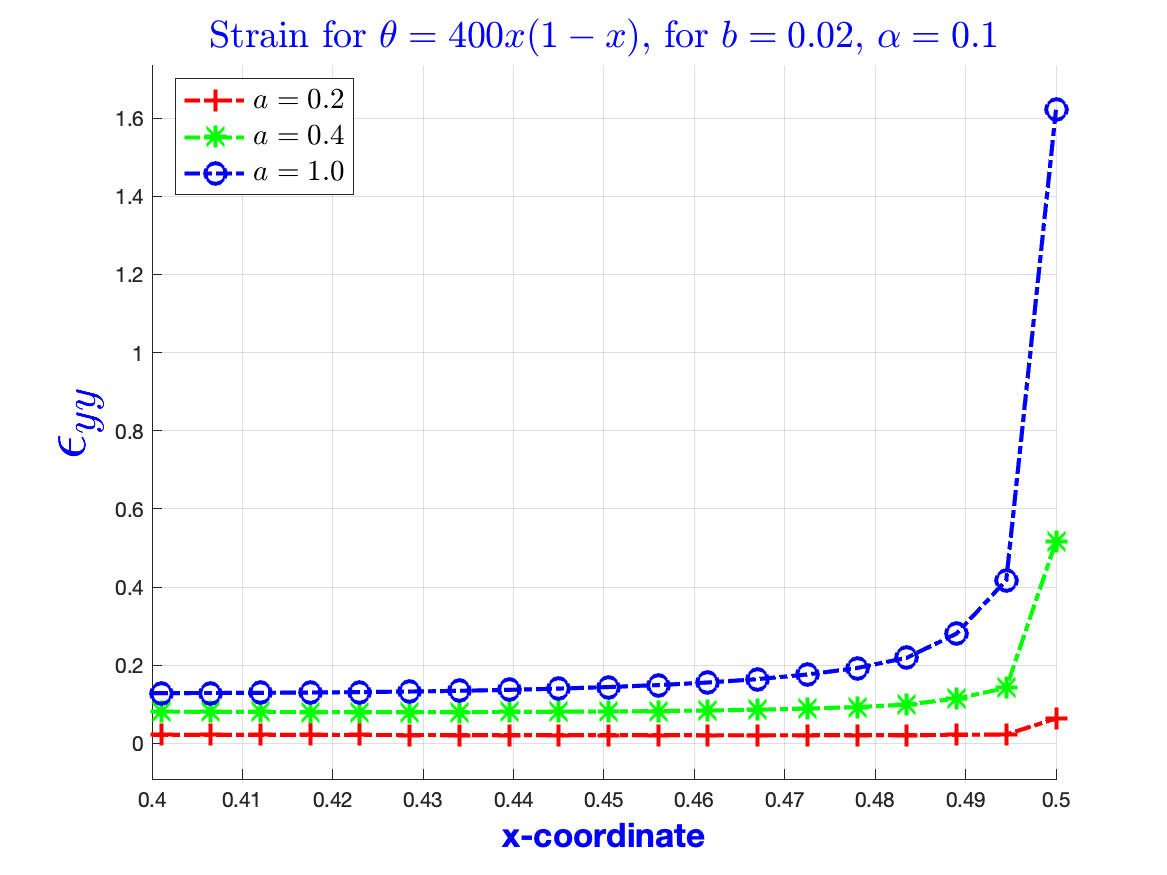}
\caption{Stress and strain for $a$ variation, $\theta=400x(1-x),b=0.02$}\label{fig:SS_m2_av_para}
\end{figure}

As shown in Figure~(\ref{fig:SS_m2_av_para}), an increase in the parameter $a$ is directly correlated with a sharp rise in stress and strain under loading. This indicates that the material stores a significantly higher amount of elastic strain energy. The rapid, uncontrolled release of this stored energy during fracture provides a powerful driving force for crack evolution, suggesting a greater propensity for {catastrophic, brittle-like failure}.

\begin{figure}[H]
\centering
\includegraphics[scale=0.2]{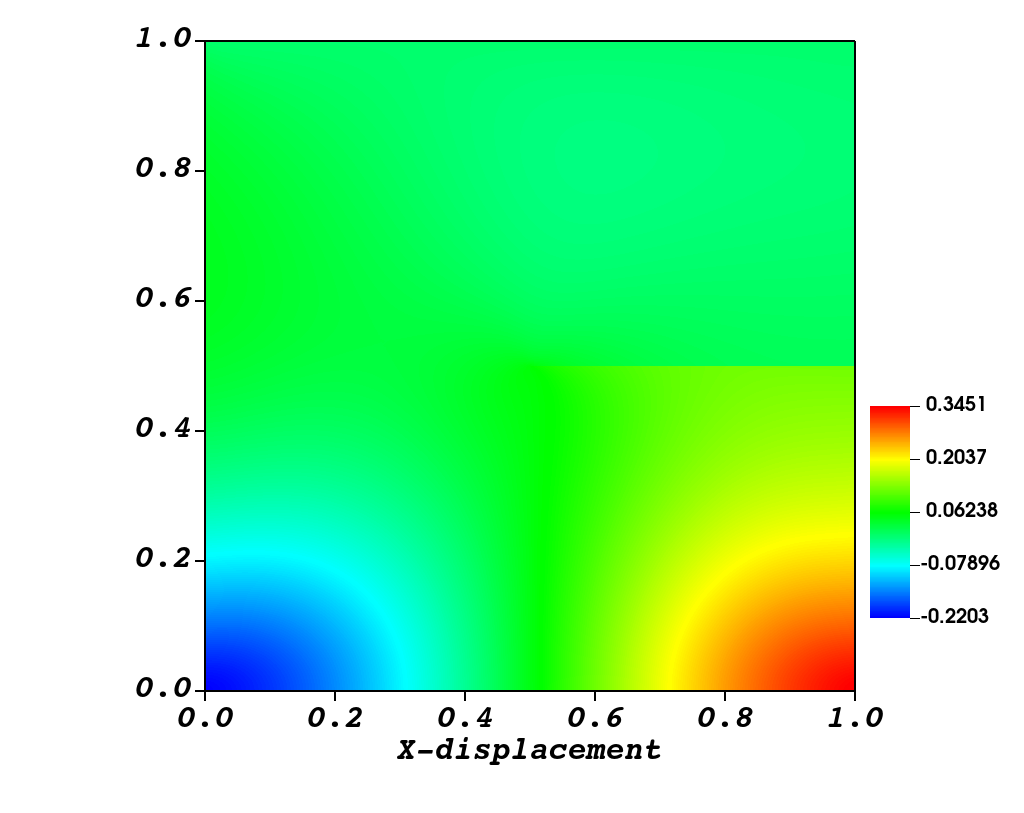}\quad
\includegraphics[scale=0.2]{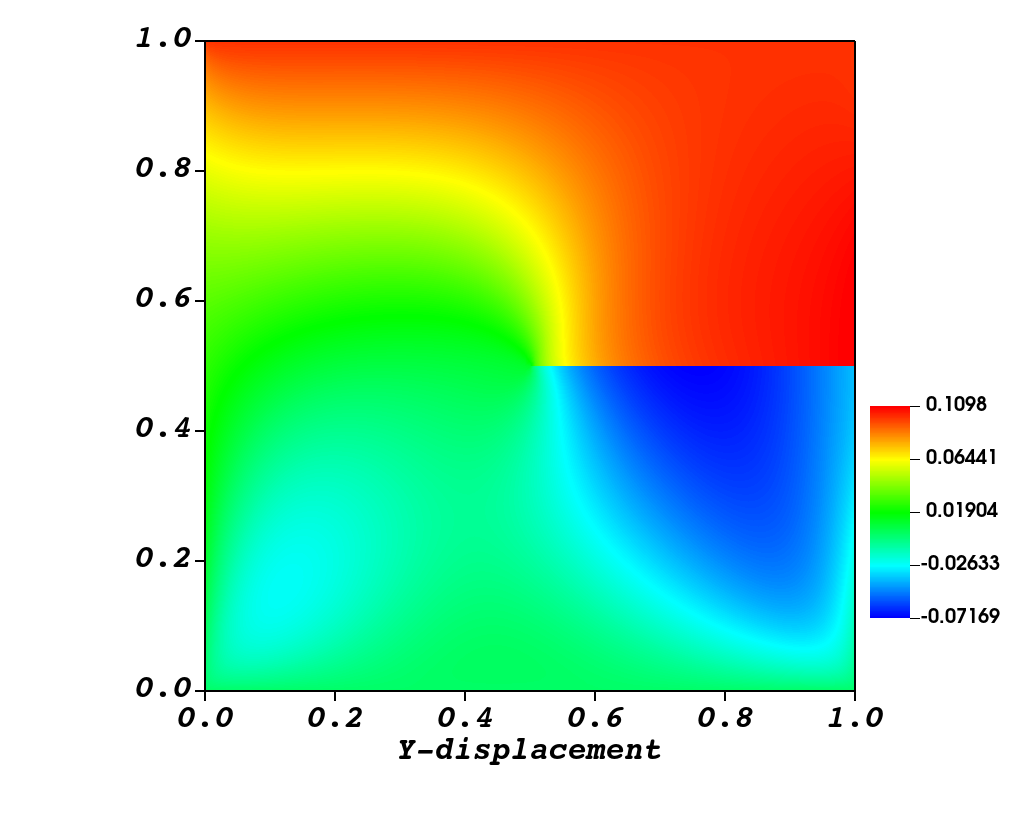}
\caption{Horizontal and vertical displacements for $\theta=400x(1-x),a=0.5,b=0.02$}\label{fig:dis_m2_para}
\end{figure}

The deformation caused by the parabolic thermal load is visualized in Figure~(\ref{fig:dis_m2_para}). An analysis of the horizontal displacement component reveals a pattern of non-uniform expansion, with a maximum value near the bottom right corner and a minimum near the bottom left. This observation is the expected outcome of the asymmetric thermal boundary condition. Furthermore, the vertical displacement field shows a peak concentration along the upper half of the right boundary. This can be interpreted as the material deflecting upwards in response to the strong compressive action generated by the lateral expansion near the base.

\begin{figure}[H]
\centering
\includegraphics[scale=0.2]{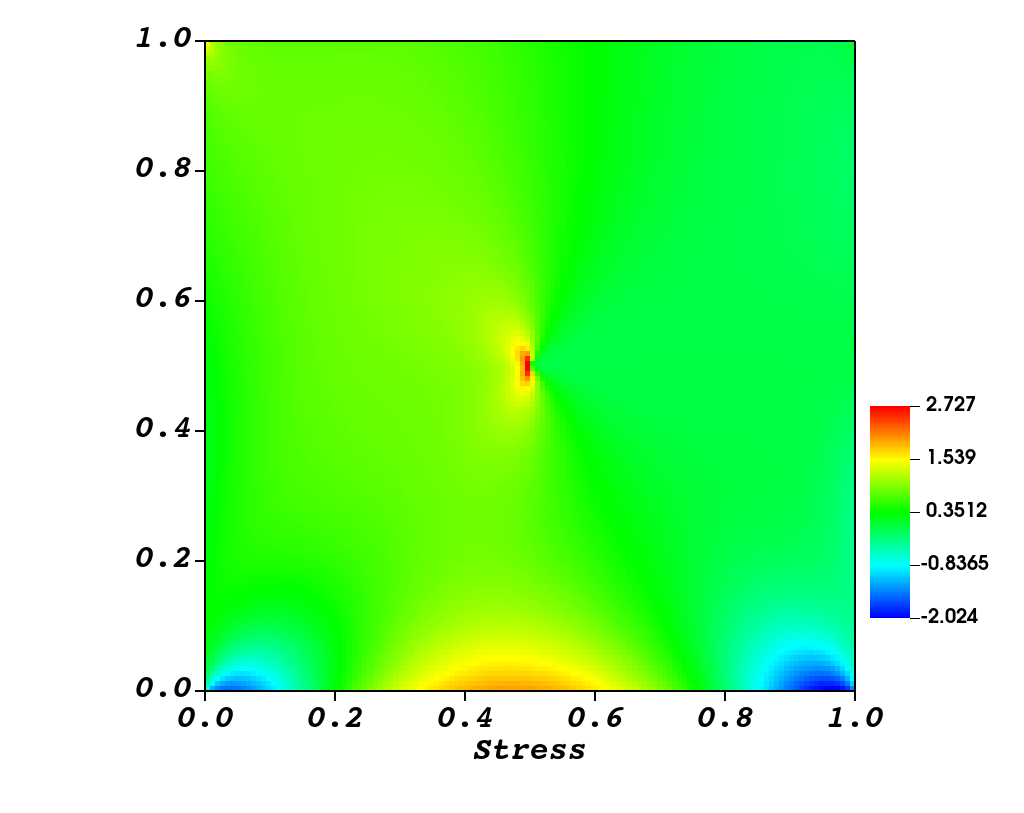}\quad
\includegraphics[scale=0.2]{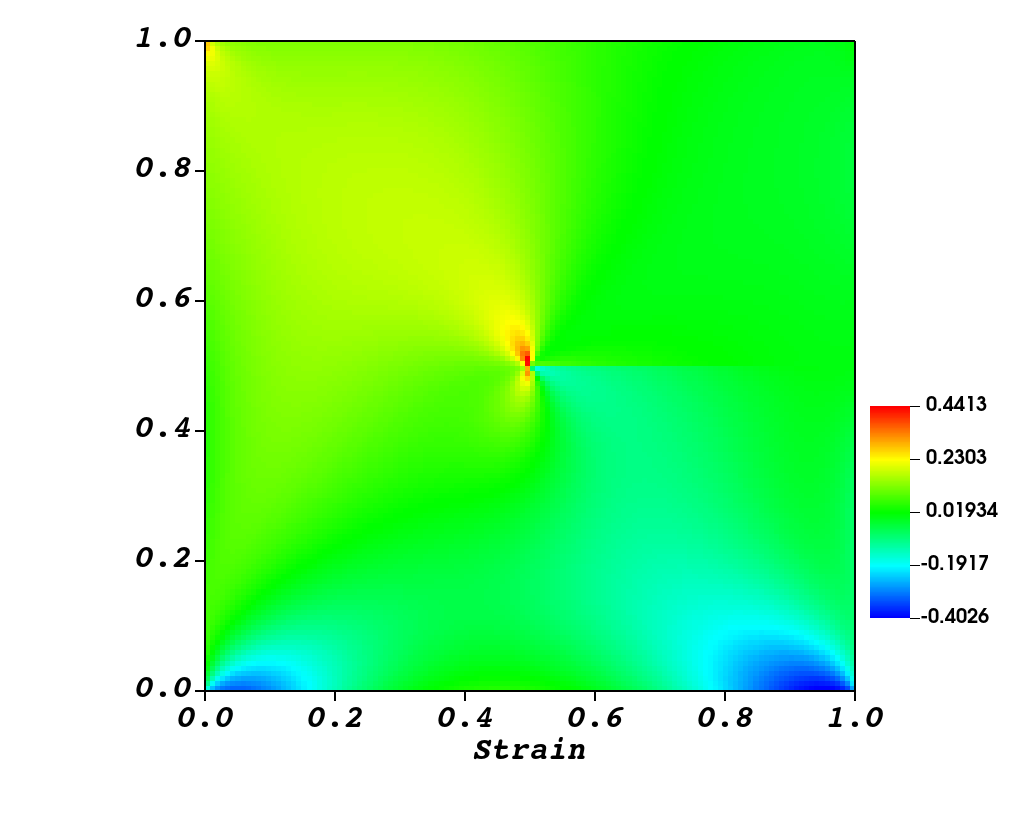}
\caption{Stress and strain for $\theta=400x(1-x),a=0.5,b=0.02$}\label{fig:SS_m2_visit_para}
\end{figure}

The stress and strain fields depicted in Figure~(\ref{fig:SS_m2_visit_para}) are dominated by a pronounced concentration at the crack tip, a characteristic feature in fracture analysis. Away from this singularity, the stress field decays smoothly and becomes relatively uniform throughout the bulk of the domain. A notable exception, however, occurs along the bottom boundary, where a localized stress amplification is evident in the central region. This secondary peak is attributed to the thermal stresses generated by the applied parabolic temperature profile. Quantitatively, the entire stress field is bounded by a maximum tensile value of $2.727$ and a minimum compressive value of $-2.024$.

\begin{figure}[H]
\centering
\includegraphics[scale=0.2]{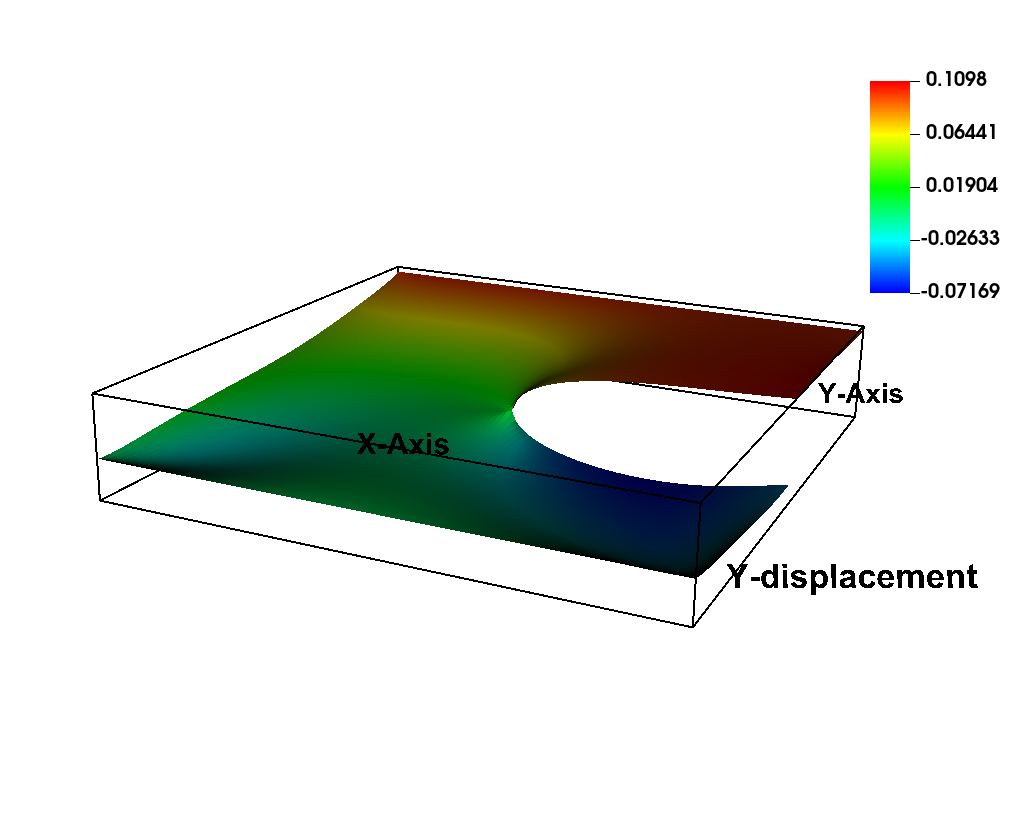} 
\centering{}\caption{Elevated $u_{y}$ for $\theta=400x(1-x),a=0.5,b=0.02$}\label{fig:y3d_m2_para}
\end{figure}

Figure~(\ref{fig:y3d_m2_para}) depicts the vertical displacement, $u_{y}$, as an elevated surface, illustrating the structural deformation caused by the parabolic temperature load on the bottom edge. The primary feature of this surface is the morphology of the crack opening, which forms a smooth, continuous profile that is approximately elliptical and centered about the crack tip. This shape represents the physical separation of the crack faces, which is driven by the non-uniform thermal expansion within the domain.

\section{Conclusion }

This study presents a finite element framework for analyzing thermoelastic
responses in transversely isotropic materials using a nonlinear strain-limiting
constitutive model. To capture the interaction between thermal and
mechanical effects, we formulate a coupled system comprising steady-state
heat equation and the quasilinear elliptic boundary value problem.
The resulting variational formulation is discretized using continuous
Galerkin methods and solved through a Picard's iteration scheme. Convergence
of the nonlinear iterations was verified by tracking the $L^{2}$
difference between successive iterates. 

Numerical experiments on edge-cracked domains, for two fiber orientations
and two different thermal boundary conditions (uniform and inhomogeneous),
reveal consistent qualitative trends. The nonlinear model preserves
the familiar crack-tip stress intensification, yet the associated
strains grow much more slowly, echoing prior observations in isotropic
thermoelasticity and purely mechanical transversely isotropic settings.
Inhomogeneous temperature fields (and their gradients) act as additional
concentrators, amplifying crack-tip fields beyond the uniform-temperature
case; this mirrors observations reported previously for isotropic
thermoelasticity.

The numerical results clearly indicate that the two primary constitutive
parameters $a$ and $b$ play contrasting roles in shaping the local
fields around the crack tip. The nonlinearity parameter $b$ acts
as a strain-limiting regularizer: as $b$ increases, the peak values
of stress and strain near the crack tip decrease consistently for
both fiber orientations and thermal boundary conditions. This behavior
suggests a toughening effect, in which the material becomes more resistant
to the concentration of stress and deformation, thereby potentially
delaying the propagation of fracture. In contrast, the parameter $a$
has an amplifying effect on the crack tip fields. Larger $a$ values
lead to higher stress magnitudes and more pronounced strain localization
at the crack tip, especially under non-uniform temperature distributions.
Together, these sensitivities underscore the importance of careful
parameter selection in the design of materials and structures where
thermomechanical fracture is a critical concern.

These findings support the relevance of strain-limiting theories in
anisotropic thermoelastic analysis and lay the groundwork for future
extensions. In particular, this framework can be leveraged to model
thermo-mechanical fracture evolution in transversely isotropic bodies
using advanced formulations such as phase-field regularization or
interface damage models.

\section{Acknoledgement}
SG would like to thank the University of Texas Rio Grande Valley for providing a Presidential Research Fellowship during his PhD studies. SMM thanks the  National Science Foundation for its financial support under Grant No. 2316905.

\bibliographystyle{plain}
\bibliography{ref_thermo}

\end{document}